\documentclass[preprint, 11pt]{elsarticle}

\usepackage[T1]{fontenc}
\usepackage[utf8]{inputenc}
\usepackage[hmargin=0.9in]{geometry}
\usepackage[babel,german=quotes]{csquotes}
\usepackage{subfig}
\usepackage{graphicx}
\usepackage[colorlinks=true,citecolor=blue,urlcolor=blue]{hyperref}
\usepackage{caption}
\usepackage{amsmath}          
\usepackage{amsthm}           
\usepackage{amssymb}          
\usepackage{bm, amsfonts}     
\usepackage{xcolor}
\usepackage{fixmath}
\usepackage{tikz}
\usepackage{bm}
\usepackage{makecell}

\newtheorem{thm}{Theorem}
\newdefinition{example}{Example}
\newproof{pf}{Proof}
\newcommand{\proofend}{\hfill $\Box$}

\newcommand {\dx} {\,{\rm d}{\mathbf x}}

\newcommand {\ds} {\,{\rm d}{\mathrm s}}

  \newcommand {\bF} {{\mathbf f}}

  \newcommand{\R}{\mathbb{R}}
  \newdefinition{rmk}{Remark}
  
  \newcommand{\pd}[2]{\frac{\partial #1}{\partial #2}}

  \newcommand{\td}[2]{\frac{\mathrm d #1}{\mathrm d #2}}

\newcommand{\beq}{\begin{equation}}
\newcommand{\eeq}{\end{equation}}

\makeatletter
\def\ps@pprintTitle{%
  \let\@oddhead\@empty
  \let\@evenhead\@empty
  \def\@oddfoot{
    \footnotesize\itshape
    \hfill\today
  }%
  \let\@evenfoot\@oddfoot}
\makeatother

\begin{document}

\begin{frontmatter} 
  \title{A new perspective on flux and slope limiting in discontinuous Galerkin methods for hyperbolic conservation laws}

\author{Dmitri Kuzmin}
\ead{kuzmin@math.uni-dortmund.de}

\address{Institute of Applied Mathematics (LS III), TU Dortmund University\\ Vogelpothsweg 87,
  D-44227 Dortmund, Germany}

\journal{Journal of Computational Physics}

\begin{abstract}
  In this work, we discuss and develop multidimensional limiting techniques for discontinuous Galerkin (DG) discretizations of scalar hyperbolic problems. To ensure that each cell average satisfies a local discrete maximum principle (DMP), we impose inequality constraints on the local Lax-Friedrichs fluxes of  a piecewise-linear ($\mathbb{P}_1$) approximation. Since the piecewise-constant ($\mathbb{P}_0$) version corresponds to a property-preserving low-order finite volume method, the validity of DMP conditions can always be enforced using slope and/or flux limiters. We show that the (currently rather uncommon) use of direct flux limiting makes it possible to construct more accurate DMP-satisfying approximations in which a weak form of slope limiting is used to prevent unbounded growth of solution gradients. Moreover, both fluxes and slopes can be limited in a manner which produces nonlinear problems with well-defined residuals even at steady state. We derive/present slope limiters based on different kinds of inequality constraints, discuss their properties and introduce new anisotropic limiters for problems that require different treatment of different space directions. At the flux limiting stage, the anisotropy of the problem at hand can be taken into account by using a customized definition of local bounds for the DMP constraints. At the slope limiting stage, we adjust the magnitude of individual directional derivatives using low-order reconstructions from cell averages to define the bounds. In this way, we avoid unnecessary limiting of well-resolved derivatives at smooth peaks and in internal/boundary layers. The properties of selected algorithms are explored in numerical studies for DG-$\mathbb{P}_1$ discretizations of two-dimensional test problems. In the context of  $hp$-adaptive DG methods, the new limiting procedures can be used in $\mathbb{P}_1$ subcells of macroelements marked as `troubled' by a smoothness indicator.
\end{abstract}
\end{frontmatter}

\section{Introduction}

Discontinuous Galerkin (DG) belong to the family of discretization techniques in which the evolution of cell averages and their values at steady state are determined by the choice of numerical fluxes. The commonly employed local Lax Friedrichs (LLF) flux approximation provably guarantees the validity of local discrete maximum principles (DMP), preservation of invariant domains, and entropy stability for piecewise-constant ($\mathbb{P}_0$) discretizations which are equivalent to cell-centered finite volume schemes. The LLF flux of a piecewise-linear ($\mathbb{P}_1$) or higher order approximation may fail to satisfy the inequality constraints that provide sufficient conditions for a DG method to possess the above properties. As a consequence, the cell averages may attain unacceptable values or violate entropy conditions.

The DMP property of a LLF flux depending on the cell averages and slopes of the DG solution in two adjacent elements can be enforced by using some form of flux and/or slope limiting. The differences between the two kinds of limiting techniques are subtle and require further explanation:
\begin{itemize}
\item A~{\it flux limiter} is an algorithm which produces a convex combination of numerical fluxes corresponding to a property-preserving low-order scheme and a high-order target discretization that may violate the DMP conditions. Classical representatives of high-resolution finite volume schemes equipped with such limiters include flux-corrected transport (FCT) algorithms \cite{shasta,zalesak79} and total variation diminishing (TVD) methods \cite{harten1,harten2} for structured grids, as well as local extremum diminishing (LED) schemes for unstructured meshes \cite{barthohl,jameson1,jameson2}. Examples of bound-preserving (BP) flux limiters for $\mathbb{P}_0$ components of
  DG solutions can be found in \cite{CH-Paper,entropyDG,Moe2017}.
  \item
    A {\it slope limiter} is an postprocessing tool which adjusts the derivatives of a piecewise-polynomial solution without changing the cell averages. The outcome of slope limiting is a convex combination of the $\mathbb{P}_0$ and $\mathbb{P}_1$ approximations in each cell. Many limiting techniques of this kind were proposed in the literature on finite volume and DG methods \cite{barthjesp,RKDG-II,giuliani,hoteit2004,kriv2007,dglim,zhang1,zhang2}. Most of them constrain the range of values that the DG solution may attain at certain control points on the boundaries of mesh elements. This limiting criterion provides safer input data for calculation of numerical fluxes but, as noticed in \cite{Moe2017}, is generally insufficient to ensure the DMP property of the cell averages that define the local bounds for the inequality constraints of the slope limiting procedure. Hence, the resulting solutions may exhibit undershoots/overshoots.
\end{itemize}

In light of the above, slope limiting can be interpreted as indirect flux
limiting via imposition of inequality constraints on the Riemann data. If the
LLF flux function is linear, a convex combination of the $\mathbb{P}_0$ and
$\mathbb{P}_1$ fluxes equals the flux corresponding to the
convex combination of the $\mathbb{P}_0$ and $\mathbb{P}_1$
 states. For linear advection in 1D,
any flux limiter is equivalent to a slope limiter and vice versa
\cite{leveque92,zalesak87}. The design of a slope limiter which ensures
the DMP property of cell averages for a general conservation law is
more involved \cite{giuliani,giuliani2020}. Many existing slope
limiters are not BP in this sense and/or may fail to preserve well-resolved
directional derivatives of the DG-$\mathbb{P}_1$ solution. In applications
to large-scale ocean flow models, the latter side effect gives rise
to spurious diapycnal mixing \cite{BurchardRennau2008,Griffies2000}
and unnecessary cancellation of all partial derivatives in boundary
elements. The use of anisotropic slope limiters
\cite{dganis,hoteit2004,berger2015} may alleviate these problems
but does not guarantee the BP property.

The poor performance of slope limiters in applications to some anisotropic
transport problems can be attributed to the fact that a steep gradient
does not produce a large flux across a surface which is (almost)
parallel to the flow direction. It is generally difficult to find
an optimal correction factor for each derivative using inequality
constraints for solution values at control points and to find the
least restrictive local bounds which guarantee the DMP property
in the absence of flux limiting. On the other hand, some control
of solution gradients is required even in DG schemes equipped
with flux limiters \cite{CH-Paper,Moe2017}. In contrast to
finite volume schemes, the difference between the unlimited
gradients of a DG-$\mathbb{P}_1$ approximation and a suitable reconstruction
from BP cell averages may become arbitrarily large. The
application of a flux limiter would prevent any
violation of local bounds but the resulting approximations may turn
out too diffusive or exhibit strong `terracing' effects, i.e., spurious
distortions within the range of values satisfying the DMP
constraints. Hence, some slope limiting is still required to
ensure consistency and stability of the DG approximation to the 
gradient but these basic requirements are far less restrictive than
inequality constraints for pointwise solution values.

The methodology that we favor in the present paper combines a
flux limiter based on inequality constraints for the cell averages
and a slope limiter based on inequality constraints for the
gradients or directional derivatives of the DG solution. The main burden
of enforcing the DMP conditions falls on the flux limiter, while
slope limiting provides a complementary correction tool that
may be applied less frequently. This unconventional
design philosophy distinguishes our new limiting strategy from mainstream
approaches which either assume the BP property of the cell averages or
enforce it using more sophisticated slope limiters than those
that we use to control the solution gradients.

The remainder of this paper is
organized as follows. In Section \ref{sec:BP}, we discretize a
generic scalar conservation law in space
using a DG-$\mathbb{P}_1$ method,
formulate the DMP constraints, and discuss the
definition of local bounds for anisotropic transport problems.
In Section \ref{sec:FL}, we review two approaches for enforcing
the corresponding inequality constraints in the process of flux
limiting. The first one is based on a localized FCT algorithm
(cf. \cite{cotter,Guermond2018,CG-BFCT}). The second type of flux
correction is based on the concept of monolithic convex limiting 
for general finite element approximations
\cite{entropyDG,convex,EG-MCL,convex2}. In Section \ref{sec:SL},
we present new slope limiters which also guarantee the DMP property
for linear and nonlinear problems. Their derivation
exploits the aforementioned relationships
between flux and slope limiting. Next, we discuss the vertex-based
version \cite{dglim} of the Barth-Jespersen slope limiter
\cite{barthjesp} for general conservation laws and propose
a more accurate slope limiting procedure. Using local
bounds defined in terms of low-order reconstructions, we constrain
the directional derivatives either via direct adjustment of their
magnitudes or via the use of penalization terms in the discretized
weak form of the governing equation.  The aspects
of entropy stabilization for nonlinear problems are briefly
discussed in Section~\ref{sec:ES}. The results of 2D numerical
experiments in Section \ref{sec:num} illustrate the typical
behavior of the presented limiting techniques in challenging
situations (anisotropic layers, propagating fronts, steady-state
conditions, nonlinear nonconvex flux functions). A grid
convergence study is performed for a test problem with a
smooth exact solution. In Section \ref{sec:conc}, we conclude
with an outlook of how the presented limiting tools can be extended
to systems and higher order space-time discretizations.

\section{General limiting framework}
\label{sec:BP}

Let $u(\mathbf{x},t)$ be a scalar conserved
quantity depending on the space location $\mathbf{x}\in
\R^d,\ d\in\{1,2,3\}$ and time moment $t\ge 0$. 
We consider the generic scalar conservation law
\beq
\pd{u}{t}+\nabla\cdot\mathbf{f}(u)=0,
\label{goveq}
\eeq
where $\mathbf{f}=(\mathsf{f}_1,\ldots,\mathsf{f}_d)$ is a (possibly
nonlinear) flux function. Let $K_i,\ i\in\{1,\ldots,E_h\}$ be a
cell of a simplex or tensor product mesh. The union $\bar\Omega_h
=\bigcup_{i=1}^{E_h}K_i$ of all mesh cells is the closure of a
bounded computational domain. The boundary $\partial K_i$ of 
the cell $K_i$ consists of {\it faces} $S_{ij}$ on which
the unit outward normal  $\mathbf{n}_{ij}$ is constant.
If $S_{ij}=\partial K_i\cap\partial K_j$ is a common
boundary of two cells, then $j\in\mathcal\{1,\ldots,E_h\}$
is the index of the adjacent cell. Boundary faces $S_{ij}\subset
\partial\Omega_h$ are numbered using
indices $j\in \{E_h+1,\ldots,\bar E_h\}$. The indices of
faces that form the boundary of $K_i$ are stored in
the set $\mathcal N_i$.

The DG-$\mathbb{P}_1$ method approximates $u\in L^2(\Omega_h)$ 
by a possibly discontinuous
piecewise-linear function $u_h\in L^2(\Omega_h)$ such that
$u_{ih}:=u_h|_{K_i}\in
\mathbb{P}_1(K_i)$ for $i=1,\ldots,E_h$.
The normal flux
$\mathbf{f}(u)\cdot\mathbf{n}_{ij}$ is approximated by
a numerical flux $H(u_{ih},u_{jh},\mathbf{n}_{ij})$. If
$S_{ij},\ j> E_h$ is a face on the boundary $\partial
\Omega_h$, then $u_{jh}=u_{ih}$ on $
\{\mathbf{x}\in S_{ij}\,:\,
\mathbf{f}'(u(\mathbf{x}))\cdot\mathbf{n}_{ij}\ge 0\}$.
For $\mathbf{f}'(u_{ih}(\mathbf{x},t))\cdot\mathbf{n}_{ij}< 0$,
the value of $u_{jh}$ on $S_{ij}$ is determined by an inflow
boundary condition. The DG-$\mathbb{P}_1$ discretization
of \eqref{goveq} on $K_i$ is given by
\beq
\label{LLF}
\int_{K_i}w_{ih}\pd{u_{ih}}{t}\dx-
\int_{K_i}\nabla w_{ih}\cdot\mathbf{f}(u_{ih})\dx
+\sum_{j\in\mathcal N_i}\int_{S_{ij}}w_{ih}
H(u_{ih},u_{jh},\mathbf{n}_{ij})\ds
=0\qquad \forall w_{ih}\in \mathbb{P}_1(K_i).
\eeq
Clearly, the accuracy and stability properties of such a method
depend on the choice of the approximate Riemann solver
$H(\cdot,\cdot,\cdot)$.
In this work, we use the local Lax-Friedrichs (LLF) flux
\beq\label{LLF_flux}
H(u_L,u_R,\mathbf{n}_{ij}) =\mathbf{n}_{ij}\cdot \frac{
  \bF(u_R)+\bF(u_L)}{2}
-\frac12\lambda_{ij}(u_L,u_R)(u_R-u_L),
\eeq
where $\lambda_{ij}$ is the maximum wave speed of
the 1D Riemann problem in the normal direction $\mathbf{n}_{ij}$, i.e.,
\beq\label{lambda}
\lambda_{ij}(u_L,u_R)=\max_{\omega\in[0,1]}|\mathbf{f}'(\omega u_R
+(1-\omega) u_L)\cdot\mathbf{n}_{ij}|.
\eeq
Using the test function $w_{ih}\equiv 1$, we find that the
evolution of the cell average
\beq U_{i0}=\frac{1}{|K_i|}\int_{K_i}u_{ih}\dx\eeq is governed
by the semi-discrete integral conservation law (cf. \cite{entropyDG,EG-MCL})
\beq\label{LLFP1}
|K_i|\td{U_{i0}}{t}+\sum_{j\in\mathcal N_i}|S_{ij}|H_{ij}^{\mathbb{P}_1}
=0,\qquad
 H_{ij}^{\mathbb{P}_1}=\frac{1}{|S_{ij}|}\int_{S_{ij}}
H(u_{ih},u_{jh},\mathbf{n}_{ij})\ds.
\eeq
The cell-centered finite volume (DG-$\mathbb{P}_0$) version of the LLF method
is defined by
\beq\label{LLFP0}
|K_i|\td{U_{i0}}{t}+\sum_{j\in\mathcal N_i}|S_{ij}|H_{ij}^{\mathbb{P}_0}
=0,\qquad
 H_{ij}^{\mathbb{P}_0}=H(U_{i0},U_{j0},\mathbf{n}_{ij}).
\eeq
This piecewise-constant space discretization preserves all
essential properties of the exact weak solution. If the system
of ordinary differential equations \eqref{LLFP0}
is discretized in time using a strong stability preserving (SSP)
Runge-Kutta method \cite{ssprev}, the resulting fully discrete
scheme is (locally) bound-preserving and satisfies a cell entropy inequality.
These well-known properties of the LLF-$\mathbb{P}_0$ method can
be easily verified using the theoretical framework developed by
Guermond and Popov \cite{Guermond2019}.

Unfortunately, the property-preserving LLF-$\mathbb{P}_0$ approximation
is first-order accurate at best. The DG-$\mathbb{P}_1$ version \eqref{LLF}
is up to second-order accurate but may converge to wrong weak solutions
and/or violate discrete maximum principles. To cure this unsatisfactory
behavior, we replace \eqref{LLF} with
\begin{align}
\int_{K_i}w_{ih}\pd{u_{ih}}{t}\dx&-
\int_{K_i}\nabla w_{ih}\cdot\mathbf{f}(u_{ih})\dx\nonumber\\
&+\sum_{j\in\mathcal N_i}\int_{S_{ij}}w_{ih}[\alpha_{ij}
  H(u_{ih}^*,u_{jh}^*,\mathbf{n}_{ij})+(1-\alpha_{ij})
H_{ij}^{\mathbb{P}_0}]\ds
=0\qquad \forall w_{ih}\in \mathbb{P}_1(K_i),\label{LLFL}
\end{align}
where $\alpha_{ij}\in[0,1]$ is an adjustable
flux correction factor such that $\alpha_{ji}=\alpha_{ij}$ and
\beq\label{slim}
u_{ih}^*(\underbrace{x_1,\ldots,x_d}_{\mathbf{x}},t)=
U_{i0}(t)+\sum_{k=1}^d\beta_{ik}U_{ik}(t)
(x_k-\bar x_{ik})
\eeq
is a limited Taylor expansion about the centroid
$\bar{\mathbf{x}}_i=\frac{1}{|K_i|}\int_{K_i}\mathbf{x}\dx$.
The multiplication of
\beq
U_{ik}:=\pd{u_{ih}}{x_k},\qquad k=1,\ldots,d
\eeq
by yet unspecified
slope correction factors $\beta_{ik}\in[0,1]$ makes it possible to
adjust the magnitude of each partial derivative if necessary.
On a non-periodic inflow boundary of $\Omega_h$, we use the external limit
$u_{jh}^*=u_{jh}$ since the value of $u_{jh},\ j>E_h$ is defined by
a Dirichlet boundary condition.

The semi-discrete evolution equation \eqref{LLFL}
corresponding to $w_{ih}\equiv 1$ blends 
\eqref{LLFP1} and \eqref{LLFP0} as follows:
\beq\label{LLFPAB}
|K_i|\td{U_{i0}}{t}+\sum_{j\in\mathcal N_i}|S_{ij}|H_{ij}^*
=0,
\eeq
\beq\label{Hstar}
H_{ij}^*=\frac{1}{|S_{ij}|}
\int_{S_{ij}}[\alpha_{ij}
  H(u_{ih}^*,u_{jh}^*,\mathbf{n}_{ij})+(1-\alpha_{ij})
H_{ij}^{\mathbb{P}_0}]\ds.
 \eeq
 Note that the limited flux $H_{ij}^*$ reduces to
 $H_{ij}^{\mathbb{P}_1}$ for $\alpha_{ij}=\beta_{ik}=\beta_{jk}=1,
 \ k=1,\ldots,d$. The flux $H_{ij}^{\mathbb{P}_0}$
 is recovered if we set
 $\alpha_{ij}=0$ or $\beta_{ik}=\beta_{jk}=0, \ k=1,\ldots,d$.
 Hence, the desired properties of the low-order LLF scheme
 \eqref{LLFP0} can always be preserved by tuning
$\alpha_{ij}$ and/or $(\beta_{ik},\beta_{jk})$.
 
 Let us discretize \eqref{LLFL} in time using an explicit SSP
 Runge-Kutta method \cite{ssprev} such that each stage has
 the structure of a forward Euler update. The fully discrete
 counterpart of
 \eqref{LLFPAB} is given by
 \beq\label{BPFE}
 U_{i0}^*=U_{i0}-\frac{\Delta t}{|K_i|}
 \sum_{j\in\mathcal N_i}|S_{ij}|H_{ij}^*,
 \eeq
 where $\Delta t>0$ is the time step. The flux $H_{ij}^*$ is
 calculated using the data from the previous time level or
 SSP-RK stage. We constrain $U_{i0}^*$ to satisfy a (local) discrete
 maximum principle of the form
 \beq\label{DMPconstr}
U_{i0}^{\min}\le  U_{i0}^*\le U_{i0}^{\max}.
 \eeq
 The choice of the bounds $U_{i0}^{\min}$ and $U_{i0}^{\max}$ must guarantee
 that the DMP constraints \eqref{DMPconstr} hold at least in the
 case of $H_{ij}^*= H_{ij}^{\mathbb{P}_0}$ for all $j\in\mathcal N_i$.
 We will define such {\it feasible} bounds in Section \ref{sec:FL}.

 Substituting $U_{i0}^*$ defined by \eqref{BPFE} 
 into \eqref{DMPconstr}, we obtain the
 equivalent inequality constraints
 \beq\label{DMPconstr2}
 \frac{|K_i|}{\Delta t}
(U_{i0}^{\min}-U_{i0})\le  -\sum_{j\in\mathcal N_i}|S_{ij}|H_{ij}^*
  \le \frac{|K_i|}{\Delta t}(U_{i0}^{\max}-U_{i0}),
 \eeq
the validity of which can be enforced via direct flux limiting
(adjustment of $\alpha_{ij}$, as discussed in Section~\ref{sec:FL})
and/or slope limiting (adjustment of $\beta_{ik},\beta_{jk}$, as
discussed in Section \ref{sec:SL}). The latter approach is commonly
preferred in DG methods (cf. \cite{giuliani,hoteit2004,dglim}). However,
instead of constraining
the linear polynomials $u_{ih}^*$
defined by \eqref{slim} in a way which guarantees that conditions
\eqref{DMPconstr2} are satisfied for $H_{ij}^*$ defined by
\eqref{Hstar}, many slope limiters impose linear inequality
constraints of the form
\beq
U_{ip}^{\min}\le u_{ih}^*(\mathbf{x}_{ip})\le U_{ip}^{\max}\eeq
on the values of $u_{ih}^*$ at certain control points $\mathbf{x}_{ip}$.
In many cases, these constraints are not equivalent to \eqref{DMPconstr2}.
If they are more restrictive, the slope limiter may fail to recognize
a well-resolved gradient. If the bounds are too wide, the resulting
cell averages $U_{ih}^*$ may violate \eqref{DMPconstr}.

In this
work, we develop both slope limiters which guarantee the validity
of DMP constraints and DG-$\mathbb{P}_1$ methods based on
the following design principles for flux and slope limiting:
\begin{itemize}
\item cell averages $U_{ih}^*$ must be kept in bounds using a flux limiter
  ($\alpha_{ij}\in[0,1]$) to enforce \eqref{DMPconstr2};
\item linear polynomials $u_{ih}^*$ may be constrained using a slope
  limiter based on inequality constraints for pointwise values
  $u_{ih}^*(\mathbf{x}_{ip})$ or partial derivatives $U_{ik}^*=
  \pd{ u_{ih}^*}{x_k}
  =\beta_{ik}\pd{ u_{ih}}{x_k},
  \ k=1,\ldots,d$.
\end{itemize}
According to this design philosophy, the purpose of slope limiting
is not to produce a flux $H_{ij}^*$ satisfying conditions
\eqref{DMPconstr2} but to
prevent unbounded growth of the partial derivatives $U_{ik}^*$. 
Failure to do so would result in an inaccurate approximation to the
cell averages $U_{i0}^*$ which represent the primary unknowns of
our method. Similarly to reconstruction-based finite volume and
DG methods, the additional Taylor degrees of freedom $U_{ik}^*$
are needed just to compute fluxes $H_{ij}^*$ that
are more accurate than $H_{ij}^{\mathbb{P}_0}$.

\section{Flux limiting}
\label{sec:FL}

Let us begin with the presentation and analysis of
flux limiters, i.e., algorithms
for calculating correction factors $\alpha_{ij}\in[0,1]$ such that
the symmetry condition $\alpha_{ij}=\alpha_{ji}$ holds and
\begin{itemize}
\item $H_{ij}^*$ defined by \eqref{Hstar} satisfies
\eqref{DMPconstr2} for a given pair of slope-limited traces
$u_{ih}^*|_{S_{ij}}$ and $u_{jh}^*|_{S_{ij}}$;
\item $H_{ji}^*=-H_{ij}^*$ satisfies similar DMP constraints
  for the cell average $U_{j0}^*$ if $1\le j\le E_h$.
\end{itemize}
Additionally, a well-designed flux limiter should produce
$\alpha_{ij}=1$ in regions where this setting does not
drive the cell averages of the DG-$\mathbb{P}_1$ solution out of bounds.

The hybrid LLF flux  \eqref{Hstar} can be written as
$H_{ij}^*=H_{ij}^{\mathbb{P}_0}-\frac{F_{ij}^*}{|S_{ij}|}$, where
$H_{ij}^{\mathbb{P}_0}$ is defined by \eqref{LLFP0} and
\beq\label{Fstar}
F_{ij}^*=\alpha_{ij}F_{ij},\quad
F_{ij}=|S_{ij}|(H_{ij}^{\mathbb{P}_0}-H_{ij}^{\mathbb{P}_1}).
\eeq
Using this representation of $H_{ij}^*$ in terms of
the low-order components $H_{ij}^{\mathbb{P}_0}$ and limited
counterparts $F_{ij}^*$ of the antidiffusive fluxes $F_{ij}$,
the DMP constraints \eqref{DMPconstr2} 
can be written as
\beq\label{DMPconstr3}
 \frac{|K_i|}{\Delta t}
(U_{i0}^{\min}-U_{i0}^{\mathbb{P}_0})\le  \sum_{j\in\mathcal N_i}F_{ij}^*
  \le \frac{|K_i|}{\Delta t}(U_{i0}^{\max}-U_{i0}^{\mathbb{P}_0}),
  \eeq
  where
 \beq\label{BPLO}
 U_{i0}^{\mathbb{P}_0}=U_{i0}-\frac{\Delta t}{|K_i|}
 \sum_{j\in\mathcal N_i}|S_{ij}|H_{ij}^{\mathbb{P}_0}.
 \eeq
 The result of flux correction depends on
 the choice of the bounds and on the
 practical algorithm for calculating correction factors
 $\alpha_{ij}\in[0,1]$ that guarantee the validity of \eqref{DMPconstr3}
 under certain time step restrictions. We elaborate on these choices 
 and analyze two simple flux limiters below.
 
\subsection{Bounds for DMP constraints}
 
 For \eqref{DMPconstr3} to be satisfied at least for
 the trivial choice $\alpha_{ij}=0\ \forall j\in\mathcal N_i$, the
 bounds $U_{i0}^{\min}$ and $U_{i0}^{\max}$ should be admissible
 in the sense that $U_{i0}^{\min}\le U_{i0}^{\mathbb{P}_0}
 \le U_{i0}^{\max}$. In particular, this will be the case
 if $U_{i0}^{\min}=u^{\min}$ and $U_{i0}^{\max}=u^{\max}$
 are some global bounds such that $ U_{i0}^{\mathbb{P}_0}
 \in  [u^{\min},u^{\max}]$. The use of global bounds
 may be appropriate, e.g., if the slope-limited linear
 polynomials $u_{ih}^*$ and $u_{jh}^*$ are constrained to
satisfy local maximum principles at certain points on $S_{ij}$
(see \cite{CH-Paper} and Section \ref{sec:SClim} below). In
this case, the target fluxes
$H_{ij}=H_{ij}^{\mathbb{P}_0}-\frac{F_{ij}}{|S_{ij}|}$
produce essentially nonoscilatory approximations and it is
sufficient to enforce the validity of a global DMP for cell
averages. For example, a flux limiter based on \eqref{DMPconstr3}
with global bounds $U_{i0}^{\min}=0$ and $U_{i0}^{\max}=1$ may be used to keep
cell-averaged concentrations in the range $[0,1]$ if
the slope limiting procedure does not
guarantee this property.

In many cases, however, a better limiting strategy is to
enforce local DMPs at the flux limiting stage. Then slope
limiting can be performed in a safe mode under less
restrictive constraints, as formulated in Section \ref{sec:DClim}.
Adopting this design philosophy, we define local bounds
\beq\label{boundsCA}
U_{i0}^{\max}=\max_{j\in \mathcal J_{i0}}U_{j0},\qquad
U_{i0}^{\min}=\min_{j\in \mathcal J_{i0}}U_{j0},
\eeq
where $\mathcal J_{i0}$ is a {\it bounding stencil} containing $i$ and (a
subset of) indices of cells $K_j$ that share a common vertex
with $K_i$. By default, the indices of all common-vertex
neighbors are included in $\mathcal J_{i0}$ but better
monotonicity preservation for anisotropic transport in
layers can be achieved by including only neighbors
belonging to the same layer in the definition of
$\mathcal J_{i0}$. The first numerical
example of Section \ref{sec:num} demonstrates the potential
benefits of using such anisotropic local bounds
for $U_{i0}^*$.

\subsection{Localized FCT limiting}

After choosing the bounds, a set of correction factors $\alpha_{ij}$
satisfying conditions \eqref{DMPconstr3} and other flux limiting criteria 
needs to be found. Introducing the bounding fluxes
\beq\label{fmaxFCT}
F_{ij}^{\max}=
\begin{cases}
  \frac{|K_i|}{\Delta t}\frac{|S_{ij}|}{|\partial K_i|}
  \max\{0,\min\{U_{i0}^{\max}-U_{i0}^{\mathbb{P}_0},
  U_{j0}^{\mathbb{P}_0}-U_{j0}^{\min}\}\}
& \mbox{if}\  1\le j\le E_h, \\[0.25cm]
\frac{|K_i|}{\Delta t}\frac{|S_{ij}|}{|\partial K_i|}
\max\{0,U_{i0}^{\max}-U_{i0}^{\mathbb{P}_0}\} & \mbox{otherwise},
\end{cases}
\eeq
\beq\label{fminFCT}
F_{ij}^{\min}=\begin{cases}
\frac{|K_i|}{\Delta t}\frac{|S_{ij}|}{|\partial K_i|}
\min\{0,\max\{U_{i0}^{\min}-U_{i0}^{\mathbb{P}_0},U_{j0}^{\mathbb{P}_0}
-U_{j0}^{\max}\}\}
& \mbox{if}\  1\le j\le E_h, \\[0.25cm]
\frac{|K_i|}{\Delta t}\frac{|S_{ij}|}{|\partial K_i|}
\min\{0,U_{i0}^{\min}- U_{i0}^{\mathbb{P}_0}\}& \mbox{otherwise},
\end{cases}
\eeq
a localized version
\cite{cotter,CG-BFCT} of Zalesak's \cite{zalesak79}
 flux-corrected transport
(FCT) algorithm defines 
\beq\label{alphaFCT}
\alpha_{ij}=\begin{cases}
\min\,\left\{1,
\frac{F_{ij}^{\max}}{F_{ij}}
\right\}
  & \mbox{if}\  F_{ij}>0,
\\[0.25cm]
 1 & \mbox{if}\  F_{ij}=0,
 \\[0.25cm]
\min\,\left\{1,
\frac{F_{ij}^{\min}}{F_{ij}}
\right\}
  & \mbox{if}\  F_{ij}<0.
\end{cases}
\eeq
Recall that $S_{ij}$ is an interior face for
$1\le j\le E_h$ and a face on the boundary $\partial\Omega_h$
for $E_h+1\le j\le \bar E_h$.

\subsection{Monolithic convex limiting}

The monolithic
convex limiting (MCL) procedure \cite{entropyDG,convex,convex2,EG-MCL}
replaces \eqref{fmaxFCT} and \eqref{fminFCT} with
\beq\label{fmaxMCL}
F_{ij}^{\max}=
\begin{cases}
\lambda_{ij}^{\mathbb{P}_0}|S_{ij}|\max\{0,
\min\{U_{i0}^{\max}-\bar U_{ij,0},\bar U_{ji,0}-U_{j0}^{\min}\}\}
& \mbox{if}\  1\le j\le E_h, \\[0.25cm]
\lambda_{ij}^{\mathbb{P}_0}|S_{ij}|\max\{0,
U_{i0}^{\max}-\bar U_{ij,0}\} & \mbox{otherwise},
\end{cases}
\eeq
\beq\label{fminMCL}
F_{ij}^{\min}=\begin{cases}
\lambda_{ij}^{\mathbb{P}_0}|S_{ij}|\min\{0,
\max\{U_{i0}^{\min}-\bar U_{ij,0},\bar U_{ji,0}
-U_{j0}^{\max}\}\}
& \mbox{if}\  1\le j\le E_h, \\[0.25cm]
\lambda_{ij}^{\mathbb{P}_0}|S_{ij}|
\min\{0,U_{i0}^{\min}- \bar U_{ij,0}\}& \mbox{otherwise},
\end{cases}
\eeq
where $\lambda_{ij}^{\mathbb{P}_0}=\lambda_{ij}
(U_{i0},U_{j0})$ is the maximum wave speed
\eqref{lambda}
of the LLF-$\mathbb{P}_0$ approximation \eqref{LLFP0} and
\beq\label{barP0}
\bar U_{ij,0}=\frac{U_{j0}+U_{i0}}{2}
-\mathbf{n}_{ij}\cdot\frac{
  \bF(U_{j0})-\bF(U_{i0})}{2\lambda_{ij}^{\mathbb{P}_0}}
\eeq
are the {\it bar states} (as originally
defined in \cite{Guermond2016})
of the property-preserving
LLF-$\mathbb{P}_0$ approximation.

Using the mean value theorem of calculus, it is easy to show that
the states $\bar U_{ij,0}$ satisfy
\cite{convex}
\beq\label{lmpbar}
\min\{U_{i0},U_{j0}\}\le \bar U_{ij,0}\le \max\{U_{i0},U_{j0}\}.
\eeq
These estimates will also hold if 
$\lambda_{ij}$ defined by \eqref{lambda}
is replaced with
an arbitrary upper bound $\lambda_{ij}^{\max}\ge\lambda_{ij}$.

The DMP properties of FCT and MCL flux limiters are summarized in
the following theorem which adapts the analysis presented in
\cite{entropyDG,convex,EG-MCL} to the setting of this paper
for the reader's convenience.

\begin{thm}[DMP property of FCT and MCL flux limiters]\label{thm1}

  Define the correction factors $\alpha_{ij}$ using formula \eqref{alphaFCT}
  with flux bounds \eqref{fmaxFCT},\eqref{fminFCT} or
  \eqref{fmaxMCL},\eqref{fminMCL} corresponding to
  \beq\label{umaxmin}
  U_{i0}^{\max}\ge\max_{j\in\mathcal N_i}U_{j0},\qquad
  U_{i0}^{\min}\le\min_{j\in\mathcal N_i}U_{j0}.
  \eeq
  Choose a time step $\Delta t$
  satisfying the CFL-like condition
  \beq\label{cfl}
\frac{\Delta t}{|K_i|}\sum_{j\in\mathcal N_i}|S_{ij}|\lambda_{ij}\le 1.
\eeq
Then update \eqref{BPFE} with limited fluxes
$H_{ij}^*$ defined by \eqref{Hstar}
satisfies the DMP constraints \eqref{DMPconstr} in the case
of periodic boundary conditions.

\end{thm} 

\begin{pf}
  The  SSP Runge-Kutta stage \eqref{BPFE} admits the following
  equivalent representations
   (cf. \cite{Guermond2016,convex})
  \begin{align}
  U_{i0}^*&=U_{i0}+\frac{\Delta t}{|K_i|}\sum_{j\in\mathcal N_i}
  (|S_{ij}|\lambda_{ij}(\bar U_{ij,0}-U_{i0})+F_{ij}^*)\label{afcform}\\
  &=U_{i0}^{\mathbb{P}_0}+\frac{\Delta t}{|K_i|}\sum_{j\in\mathcal N_i}
  F_{ij}^*\label{fctform}\\
  &=U_{i0}+\frac{\Delta t}{|K_i|}\sum_{j\in\mathcal N_i}
  |S_{ij}|\lambda_{ij}(\bar U_{ij,0}^*-U_{i0}),\label{fcbarform}
  \end{align}
  where $U_{i0}^{\mathbb{P}_0}=U_{i0}+\frac{\Delta t}{|K_i|}\sum_{j\in\mathcal N_i}
  |S_{ij}|\lambda_{ij}(\bar U_{ij,0}-U_{i0})
  $ is the result of a forward Euler step for
  \eqref{LLFP0} and
  \beq\label{ubarsstarij}
  \bar U_{ij,0}^*=\bar U_{ij,0}+\frac{F_{ij}^*}{|S_{ij}|\lambda_{ij}}=
  \bar U_{ij,0}+\frac{\alpha_{ij}F_{ij}}{|S_{ij}|\lambda_{ij}}.
  \eeq
  are flux-corrected counterparts of the $\mathbb{P}_0$ bar
  states $\bar U_{ij,0}$ defined by \eqref{barP0}. Invoking
  \eqref{lmpbar} and \eqref{umaxmin}, we find that
  $\bar U_{ij,0}\in[U_{i0}^{\min},U_{i0}^{\max}]$. If the time
  step $\Delta t$ satisfies \eqref{cfl} then
   \beq\label{barstar}
  U_{i0}^*
=\left(1-\frac{\Delta t}{|K_i|}\sum_{j\in\mathcal N_i}
  |S_{ij}|\lambda_{ij}\right)U_{i0}+\frac{\Delta t}{|K_i|}\sum_{j\in\mathcal N_i}
  |S_{ij}|\lambda_{ij}\bar U_{ij,0}^{*}
  \eeq
  is a convex combination of $U_{i0}$ and $\bar U_{ij,0}^{*}$.
  Therefore, update \eqref{barstar} is bound-preserving
  if the choice of $\alpha_{ij}$ in \eqref{ubarsstarij}
 guarantees that $U_{ij,0}^{*}\in[U_{i0}^{\min},U_{i0}^{\max}]$ for
 $j\in\mathcal N_i$. In particular, this sufficient
 condition is satisfied for the trivial choice $\alpha_{ij}=0$.
 It follows that $U_{i0}^{\mathbb{P}_0}\in[U_{i0}^{\min},U_{i0}^{\max}]$.

 If $\alpha_{ij}$ are defined by
 \eqref{alphaFCT} with the MCL 
 bounds \eqref{fmaxMCL},\eqref{fminMCL} then the BP property of
 $U_{ij,0}^{*}$ follows from the estimates
 $F_{ij}^{\min}\le \alpha_{ij}F_{ij}\le F_{ij}^{\max}$
 by definition of $F_{ij}^{\min}$ and $F_{ij}^{\max}$,
 as originally shown in \cite{convex}. For
 $\alpha_{ij}$ defined by
 \eqref{alphaFCT} with the FCT 
 bounds \eqref{fmaxFCT},\eqref{fminFCT}, the limited
 fluxes satisfy
 \beq\label{fctlocest}
\frac{|K_i|}{\Delta t}\frac{|S_{ij}|}{|\partial K_i|}
(U_{i0}^{\min}-U_{i0}^{\mathbb{P}_0})\le F_{ij}^*=
\alpha_{ij}F_{ij}
 \le\frac{|K_i|}{\Delta t}\frac{|S_{ij}|}{|\partial K_i|}
(U_{i0}^{\max}-U_{i0}^{\mathbb{P}_0})
 \eeq
 with $\alpha_{ij}\in[0,1]$
 provided that $U_{i0}^{\mathbb{P}_0}\in[U_{i0}^{\min},U_{i0}^{\max}]$.
 As shown above, the low-order LLF approximation $U_{i0}^{\mathbb{P}_0}$
 is in bounds for time steps satisfying \eqref{cfl}.
 Substituting estimates \eqref{fctlocest} into \eqref{fctform},
 we find that the assertion of the theorem is true
 for the localized FCT limiter as well. \proofend
\end{pf}

\begin{rmk}
Condition \eqref{umaxmin} may be violated if $\mathcal N_i$ is not
a subset of the set $\mathcal J_{i0}$ that we used to define the
local bounds in \eqref{boundsCA}, i.e., if not all common face
neighbors are included to use the option of anisotropic flux
limiting. In this case, our definitions \eqref{fmaxFCT},\eqref{fminFCT}
and \eqref{fmaxMCL},\eqref{fminMCL} of the bounding fluxes automatically
extend the bounds to include the admissible states $U_{i0}^{\mathbb{P}_0}$
and $\bar U_{ij,0}$, respectively.
  \end{rmk}

\begin{rmk}
A linear flux function $\mathbf{v}u$ depending on a spatially variable
velocity field $\mathbf{v}=\mathbf{v}(\mathbf{x})$ is not of the
form $\mathbf{f}(u)$. As shown in \cite{EG-MCL}, definition \eqref{barP0}
of the bar states $\bar U_{ij,0}$ should be replaced with
$$
\bar U_{ij,0}=\frac{U_{j0}+U_{i0}}{2}
-\left(\frac{\int_{S_{ij}}\mathbf{v}\cdot\mathbf{n}_{ij}\ds}{
       \int_{S_{ij}}|\mathbf{v}\cdot\mathbf{n}_{ij}|\ds}\right)
  \frac{U_{j0}-U_{i0}}{2}
$$
in the case of the linear advection equation. 
Moreover, the additional `reactive' term
$$-\frac{\Delta t}{|K_i|}U_{i0}
\sum_{j\in\mathcal N_i}\int_{S_{ij}}\mathbf{v}\cdot\mathbf{n}_{ij}\ds
$$
  appears on the right-hand sides of \eqref{afcform}--\eqref{fcbarform}.
  If the vector field $\mathbf{v}$ is not divergence-free, only
  positivity preservation can be guaranteed for the
  exact and numerical solution.
\end{rmk}

  A detailed description of FCT and MCL flux limiters for the linear
  advection equation (including DMP analysis for non-periodic
  boundary conditions) can be found in \cite{EG-MCL}, where
  such limiters were used to constrain the cell averages of a $\mathbb{P}_1
  \oplus\mathbb{P}_0$ enriched Galerkin approximation. The two limiting
  approaches produce very similar results for time-dependent hyperbolic
  problems but the MCL version is better suited for calculating stationary
  solutions because the flux bounds \eqref{fmaxMCL},\eqref{fminMCL} are
  independent of the (pseudo-)time step $\Delta t$, and so is the fixed
  point solution to which \eqref{BPFE} converges. For that reason,
  we perform flux limiting using the MCL version of \eqref{alphaFCT}
  in the numerical experiments of Section \ref{sec:num}.

\section{Slope limiting}
\label{sec:SL}

The traditional purpose of slope limiting in DG methods is to
enforce inequality constraints for pointwise values of numerical
solutions {\it assuming} that the cell averages are
bound-preserving. This assumption is generally not true
but slope limiting preconstrains the numerical fluxes in
a way which makes violations of \eqref{DMPconstr} less
likely or even provably impossible. In this section, we
discuss three ways to construct slope limiters. The first
approach guarantees the validity of \eqref{DMPconstr2}. The
limiting formula for the linear advection equation is
relatively simple and equivalent to \eqref{alphaFCT}
in the 1D case. The one for nonlinear problems is derived
using linear sufficient conditions. The
second slope limiter ensures that the value of the
linear polynomial $u_{ih}^*$ at each vertex of $K_i$
is bounded by the maximum and minimum of $U_{j0}$
  in cells $K_j$
  containing the vertex \cite{dglim}. This limiter is generally
not BP but produces nonoscillatory solutions even in the
absence of flux limiting. The third approach that we
consider is slope limiting based on inequality constraints
for the partial derivatives $U_{ik}^*=\pd{ u_{ih}^*}{x_k}$
rather than pointwise values of~$u_{ih}^*$. Delegating
enforcement of the DMP constraints \eqref{DMPconstr}
to the flux limiters presented in
Section~\ref{sec:FL}, this algorithm keeps $U_{ik}^*,\ k=1,\ldots,d$
in bounds that depend on low-order reconstructions from
cell averages. The first two limiters are {\it isotropic},
i.e., they use $\beta_{i1}=\cdots=\beta_{id}=\beta_i\in[0,1]$
in \eqref{slim}. The third slope limiter is {\it anisotropic}
in the sense that individually chosen correction factors
$\beta_{ik}$ are applied to $U_{ik},\ k=1,\ldots,d$.
A monolithic version of this limiter
penalizes violations of the inequality constraints for
$U_{ik}^*$ without manipulating the partial derivatives
$U_{ik}$ of
the DG-$\mathbb{P}_1$ solution directly. 

\subsection{Isotropic slope limiting under flux constraints}

We begin with the derivation of a slope limiter based on
the flux constraints \eqref{DMPconstr2} for a linear
conservation law.
If the flux limiter is deactivated by setting $\alpha_{ij}=1\
\forall j\in\mathcal N_i$, and an isotropic slope limiting 
strategy is adopted, the limited flux $H_{ij}^*,\ 1\le j\le E_h$
defined by \eqref{slim} and \eqref{Hstar} reduces to
\beq\label{SLfux}
H_{ij}^*=\frac{1}{|S_{ij}|}
\int_{S_{ij}}H\left(U_{i0}+\beta_i(u_{ih}-U_{i0}),
U_{j0}+\beta_j(u_{jh}-U_{j0}),\mathbf{n}_{ij}\right)\ds,
\eeq
where $u_{ih}\in \mathbb{P}_1(K_i)$ is the unconstrained DG-$\mathbb{P}_1$
solution, $U_{i0}$ is its average in $K_i$ and $\beta_i\in[0,1]$
is the isotropic slope limiting factor for the constant gradient of
$u_{ih}$.

In the case of a linear advection equation with 
$\mathbf{f}(u)=\mathbf{v}u$, the LLF flux becomes the upwind flux
\beq
H(u_L,u_R,\mathbf{n})=
\begin{cases}
  (\mathbf{v}\cdot\mathbf{n})u_R & \mbox{if}\ \mathbf{v}\cdot\mathbf{n}<0,\\
  0 & \mbox{if}\ \mathbf{v}\cdot\mathbf{n}=0,\\
  (\mathbf{v}\cdot\mathbf{n})u_L & \mbox{if}\ \mathbf{v}\cdot\mathbf{n}> 0
 \end{cases}
\eeq
and, therefore, the corresponding limited flux $H_{ij}^*$ can be written as
\beq\label{HstarLA}
H_{ij}^*=
H_{ij}^{\mathbb{P}_0}+\frac{\beta_i}{|S_{ij}|}
\int_{S_{ij}\cap\partial K_+}(u_{ih}-U_{i0})\mathbf{v}\cdot\mathbf{n}_{ij}\ds
+\frac{\beta_j}{|S_{ij}|}
\int_{S_{ij}\cap\partial K_-}(u_{jh}-U_{j0})\mathbf{v}\cdot\mathbf{n}_{ij}\ds,
\eeq
where $\partial K_+=\{\mathbf{x}\in\partial K\,:\:\mathbf{v}\cdot
\mathbf{n}>0\}$ and $\partial K_-=\{\mathbf{x}\in \partial K\,:\:\mathbf{v}\cdot
\mathbf{n}<0\}$ are the outlet and inlet of $\partial K_i$.

It follows that the validity of the inequality constraints
\eqref{DMPconstr3} for $F_{ij}^*:=F_{ij}$, where
\beq\label{flimSL1int}
F_{ij}=\beta_i
\int_{S_{ij}\cap \partial K_+}(U_{i0}-u_{ih})\mathbf{v}\cdot\mathbf{n}_{ij}\ds
+\beta_j
\int_{S_{ij}\cap\partial K_-}(U_{j0}-u_{jh})\mathbf{v}\cdot\mathbf{n}_{ij}\ds
\eeq
can be enforced by using the slope correction factors 
\beq\label{betaSL1}
\beta_i=\min_{j\in\mathcal N_i^+}
\begin{cases}
  \alpha_{ij}^+ & \mbox{if}\
  \int_{S_{ij}\cap\partial K_+}(U_{i0}-u_{ih})\mathbf{v}\cdot\mathbf{n}_{ij}\ds >0,\\ 
  1 & \mbox{if}\
  \int_{S_{ij}\cap\partial K_+}(U_{i0}-u_{ih})\mathbf{v}\cdot\mathbf{n}_{ij}\ds =0,\\ 
  \alpha_{ij}^- & \mbox{if}\
  \int_{S_{ij}\cap\partial K_+}(U_{i0}-u_{ih})\mathbf{v}\cdot\mathbf{n}_{ij}\ds <0,
  \end{cases}
\eeq
where $\mathcal N_i^+$ is the set of faces $S_{ij}$ such
that $|S_{ij}\cap \partial K_+|>0$ for
$\partial K_+=\{\mathbf{x}\in\partial K\,:\:\mathbf{v}\cdot
\mathbf{n}>0\}$ and
\beq\label{alphaSL1}
\alpha_{ij}^+=\min\,\left\{1,
\frac{F_{ij}^{\max}}{F_{ij}^+}\right\},\qquad
\alpha_{ij}^-=\min\,\left\{1,
\frac{F_{ij}^{\min}}{F_{ij}^-}\right\}
\eeq
are the face-based correction factors that formula
\eqref{alphaFCT} would produce for the antidiffusive fluxes
\begin{align}\label{fmaxSL}
F_{ij}^+&=\max\left\{0,
\int_{S_{ij}\cap\partial K_+}(U_{i0}-u_{ih})\mathbf{v}\cdot\mathbf{n}_{ij}\ds\right\}
+\max\left\{0,\int_{S_{ij}\cap\partial K_-}(
U_{j0}-u_{jh})\mathbf{v}\cdot\mathbf{n}_{ij}\ds\right\}, \\
F_{ij}^-&=\min\left\{0,
\int_{S_{ij}\cap\partial K_+}(U_{i0}-u_{ih})
\mathbf{v}\cdot\mathbf{n}_{ij}\ds\right\}
+\min\left\{0,\int_{S_{ij}\cap\partial K_-}(
U_{j0}-u_{jh})\mathbf{v}\cdot\mathbf{n}_{ij}\ds\right\}.\label{fminSL}
\end{align}
The structure of this new isotropic slope limiter is similar to
that of Zalesak's FCT flux limiter \cite{zalesak79}.

\begin{rmk}
  If $S_{ij}\subset\partial\Omega_h$, which is the case for
  $E_h+1\le j\le\bar E_h$, and a Dirichlet boundary data
  $u_{jh}$ is prescribed at the inlet $S_{ij}\cap
\partial K_-=\{\mathbf{x}\in S_{ij}\,:\:\mathbf{v}\cdot
\mathbf{n}<0\}$,  then we have
\beq\label{flimSL1bdr}
F_{ij}=\beta_i
\int_{S_{ij}\cap\partial K_+}(U_{i0}-u_{ih})\mathbf{v}\cdot\mathbf{n}_{ij}\ds
\eeq
and the boundary integrals over $S_{ij}\cap\partial K_-$
should be omitted in definitions \eqref{fmaxSL}, \eqref{fminSL}.
\end{rmk}

\begin{thm}[DMP property of the new slope limiter with linear flux constraints]\label{thm2}
  Let slope limiting be performed using formula \eqref{betaSL1}
  under assumptions of Theorem \ref{thm1}. Then
  update~\eqref{BPFE} with $H_{ij}^*$ defined by \eqref{HstarLA}
  satisfies the DMP constraints \eqref{DMPconstr}.
\end{thm} 

\begin{pf}
  Using \eqref{flimSL1int}--\eqref{fminSL},
  we obtain the estimates $F_i^{\min}\le\alpha_{ij}^-F_{ij}^-\le
  F_{ij}\le\alpha_{ij}^+F_{ij}^+\le F_{ij}^{\max}$. The
  validity of the discrete maximum principle
  follows as in the proof
  of Theorem \ref{thm1}. 
   \qquad\proofend
\end{pf}

\begin{rmk}
For the 1D linear advection equation $\pd{u}{t}+v\pd{u}{x}=0$ with constant
velocity $v\in\mathbb{R}\backslash\{0\}$, the set $\mathcal N_i^+$
consists of a single index $j$. In this case, formula
\eqref{betaSL1} produces $\beta_i=\alpha_{ij}$, where $\alpha_{ij}$
is given by \eqref{alphaFCT} with the upwind-sided antidiffusive
flux $F_{ij}=|v|(U_{i0}-u_{ih}|_{S_{ij}})=H_{ij}^{\mathbb{P}_0}-H_{ij}^{\mathbb{P}_1}$.
Hence, the slope limiter defined by \eqref{betaSL1} is equivalent
to the generic flux limiter defined by \eqref{alphaFCT}.
\end{rmk}

In the case of a nonlinear conservation law, the limited flux
\eqref{SLfux} depends on $\beta_i$ and $\beta_j$ in a
nonlinear manner. However, linearized sandwich estimates
of the form (cf. \cite{FCTsync})
\beq\label{DMPconstr5}
F_{ij}^{\min}\le \alpha_{ij}^-F^-_{ij}\le 
F_{ij}(\beta_i,\beta_j)\le \alpha_{ij}^+F_{ij}^+
\le F_{ij}^{\max}\qquad\forall \beta_i,\beta_j\le \min\{\alpha_{ij}^+,\alpha_{ij}^-\}
\eeq
with $\alpha_{ij}^\pm$ defined by \eqref{alphaSL1}
can be used to enforce \eqref{DMPconstr} again. Introducing
\beq
P_{ij}^+=\int_{S_{ij}}\max\{0,U_{i0}-u_{ih}\}\ds,
\qquad P_{ij}^-=\int_{S_{ij}}\min\{0,U_{i0}-u_{ih}\}\ds,
\eeq
we define the correction factor $\beta_i$ of the slope limiter with
nonlinear flux constraints as follows:
\beq\label{betaSL1nonlin}
\beta_i=\min_{j\in\mathcal N_i}
\begin{cases}
  \min\{\alpha_{ij}^+,\alpha_{ij}^-\} & \mbox{if}\ P_{ij}^+>0,\ P_{ij}^-<0,\\
  \alpha_{ij}^+ & \mbox{if}\ P_{ij}^+>0,\ P_{ij}^-=0,\\
  \alpha_{ij}^- & \mbox{if}\ P_{ij}^+=0,\ P_{ij}^-<0,\\
  1 & \mbox{if}\ P_{ij}^+=0,\ P_{ij}^-=0.
  \end{cases}
\eeq
To derive the bounding fluxes $F_{ij}^\pm$ for \eqref{DMPconstr5}, we
need to estimate $F_{ij}(\beta_i,\beta_j)$. Let us define the
general LLF flux
\eqref{LLF_flux} for \eqref{SLfux} using
an upper bound  $\lambda_{ij}^{\max}$ for the maximum
wave speed such that
\beq
\lambda_{ij}^{\max}\ge \lambda_{ij}(U_{i0}+\beta_i(u_{ih}-U_{i0}),
U_{j0}+\beta_j(u_{jh}-U_{j0}))\quad\forall \beta_i,\beta_j\in[0,1].
\eeq
Then the slope-limited antidiffusive flux is given by
$F_{ij}(\beta_i,\beta_j)=\int_{S_{ij}}G(\beta_i,\beta_j,\mathbf{n}_{ij})\ds$,
where
\begin{align}
G(\beta_i,\beta_j,\mathbf{n}_{ij})&=\mathbf{n}_{ij}\cdot
\frac{\mathbf{f}(U_{j0})-
  \mathbf{f}(U_{j0}+\beta_j(u_{jh}-U_{j0}))}{2}\nonumber\\&+
\mathbf{n}_{ij}\cdot
\frac{\mathbf{f}(U_{i0})-
  \mathbf{f}(U_{i0}+\beta_i(u_{ih}-U_{i0}))}{2}\nonumber\\
&-\frac{\lambda_{ij}^{\max}}{2}[
\beta_j(U_{j0}-u_{jh})-\beta_i(U_{i0}-u_{ih})
  ].
\end{align}
By the mean value theorem, there exist intermediate states
$u_R$ and $u_L$ such that
\begin{align}
\mathbf{n}_{ij}\cdot[\mathbf{f}(U_{j0})-
    \mathbf{f}(U_{j0}+\beta_j(u_{jh}-U_{j0}))]
&=\beta_j\mathbf{n}_{ij}\cdot
\mathbf{f}'(u_R)(u_{jh}-U_{j0}),\\
\mathbf{n}_{ij}\cdot[\mathbf{f}(U_{i0})-
    \mathbf{f}(U_{i0}+\beta_i(u_{ih}-U_{i0}))]
& =\beta_i\mathbf{n}_{ij}\cdot
\mathbf{f}'(u_L)(u_{ih}-U_{i0}),
\end{align}
where
$|\mathbf{n}_{ij}\cdot\mathbf{f}'(u_R))|\le\lambda_{ij}^{\max}$ and
$|\mathbf{n}_{ij}\cdot\mathbf{f}'(u_L))|\le\lambda_{ij}^{\max}$ 
by definition of $\lambda_{ij}^{\max}$ as an upper bound
for the maximum wave speed. It follows that
conditions \eqref{DMPconstr5} are satisfied for the bounding fluxes
\begin{align}
  F_{ij}^+&=\lambda_{ij}^{\max}\int_{S_{ij}}(
 \max\{0,U_{i0}-u_{ih}\}-\min\{0,U_{j0}-u_{jh}\})\ds=
\lambda_{ij}^{\max}(P_{ij}^+-P_{ji}^-),\\
F_{ij}^-&=\lambda_{ij}^{\max}\int_{S_{ij}}(
\min\{0,U_{i0}-u_{ih}\}-\max\{0,U_{j0}-u_{jh}\})\ds
=\lambda_{ij}^{\max}(P_{ij}^--P_{ji}^+).
\end{align}

\begin{thm}[DMP property of the new slope limiter with nonlinear flux constraints]\label{thm3}
  Let slope limiting be performed using formula \eqref{betaSL1nonlin}
  under assumptions of Theorem \ref{thm1}. Then update~\eqref{BPFE}
  with $H_{ij}^*$ defined by \eqref{Hstar} satisfies the DMP
  constraints \eqref{DMPconstr} for any choice of $\alpha_{ij}\in[0,1]$.
\end{thm} 

\begin{pf}
 The proof is similar to that of
  of Theorem \ref{thm2}. 
   \qquad\proofend
\end{pf}

\begin{rmk}
In principle, the nonlinear flux constraints
$F_{ij}^{\min}\le F_{ij}(\beta_i,\beta_j)\le F_{ij}^{\max}$
can be enforced using optimization-based 
methods to determine the values of $\beta_i$ and $\beta_j$.
However, the high cost of solving nonlinear
inequality-constrained optimization problems makes
this approach impractical since the same DMP constraints
can be easily enforced using the flux limiters presented
in Section \ref{sec:FL}.
\end{rmk}

In summary, slope limiting under flux constraints can guarantee
the DMP property of cell averages in the nonlinear
case as well. Unfortunately, the resulting
schemes are either costly or based on pessimistic estimates.
The slope limiters to be presented next are based on linear
inequality constraints. These limiters do not control the
cell averages directly but provide sufficiently accurate
input for subsequent flux limiting, which is our preferred
approach to keeping the cell averages in bounds.

\subsection{Isotropic slope limiting under solution constraints}
\label{sec:SClim}

Instead of enforcing the flux constraints \eqref{DMPconstr2}, the
slopes of the DG solution can be adjusted to ensure that the value
of the linear polynomial $u_{ih}^*=U_{i0}+\beta_i(u_{ih}-U_{i0})$
at any point $\mathbf{x}\in K_i$
will be bounded by the maximum and minimum of cell averages in
surrounding cells \cite{dglim}. Since $u_{ih}^*$ attains
its extrema at the vertices $\mathbf{x}_{i1},\ldots,\mathbf{x}_{iN}$
of $K_i$, it is sufficient to ensure that the pointwise values
$u_{ih}(\mathbf{x}_{ip}),\ p=1,\ldots,N$ are in bounds. The
vertex-based DG version \cite{dglim} of the Barth-Jespersen
slope limiter \cite{barthjesp} for finite volume schemes
is designed to enforce inequality constraints of the form
\beq\label{vblmp}
U_{ip}^{\min}
\le u_{ih}^*(\mathbf{x}_{ip})=U_{i0}+\beta_i(u_{ih}(\mathbf{x}_{ip})-U_{i0})
\le  U_{ip}^{\max},\qquad p=1,\ldots,N
\eeq
using the correction factor
\beq\label{betadef}
\beta_i=
\min_{1\le p \le N}\begin{cases}
  \min\left\{1,\frac{U_{ip}^{\max}-U_{i0}}{u_{ih}(\mathbf{x}_{ip})-U_{i0}}
    \right\} & \mbox{if}\  u_{ih}(\mathbf{x}_{ip})>U_{i0},\\
           1 &\mbox{if}\ u_{ih}(\mathbf{x}_{ip})=U_{i0},\\
           \min\left\{1,\frac{U_{ip}^{\min}-U_{i0}}{u_{ih}(\mathbf{x}_{ip})-U_{i0}}
    \right\} & \mbox{if}\  u_{ih}(\mathbf{x}_{ip})<U_{i0}.
    \end{cases}
\eeq
Remarkably, this formula has the same structure as definition
\eqref{alphaFCT} of the flux correction factors $\alpha_{ij}$.

By default, the local bounds $U_{ip}^{\min}$ and $U_{ip}^{\max}$ of the
vertex-based (VB) slope limiter are defined by
\beq\label{boundsVB}
U_{ip}^{\max}=\max_{j\in\mathcal E_{ip}}U_{j0},\qquad
U_{ip}^{\min}=\min_{j\in\mathcal E_{ip}}U_{j0},
\eeq
where $\mathcal E_{ip}$ is the integer set containing the indices of cells that meet
at the vertex $\mathbf{x}_{ip},\ p\in\{1,\ldots,N\}$. 

As demonstrated, e.g., in \cite{Beisiegel}, the isotropic VB slope limiter
defined by \eqref{vblmp} and \eqref{betadef} belongs to the most accurate
limiting techniques for DG schemes. It is easy to implement and typically
produces solutions that are free of undershoots and overshoots. However,
it is not provably bound-preserving if no flux limiting is performed
to enforce the DMP property \eqref{DMPconstr} of the cell averages
that define the local bounds \eqref{boundsVB}.
Indeed, the validity of the flux constraints \eqref{DMPconstr2}
for $H_{ij}^*$ defined by \eqref{Hstar} with $\alpha_{ij}=1$ for all
$j\in\mathcal N_i$ does not follow from \eqref{vblmp}, although
a violation of these constraints is unlikely in practice. In
fact, the isotropic VB limiter tends to overconstrain $u_{ih}^*$.
This tendency manifests itself, e.g., in the unnecessary cancellation
of the solution gradients in boundary elements, where the default local
bounds \eqref{boundsVB} of the inequality
constraints~\eqref{betadef} are too restrictive for $\mathbf{x}_{ip}\in\partial
\Omega_h$. A possible remedy to this drawback is a customized
definition of $U_{ip}^{\max}$ and  $U_{ip}^{\min}$ for boundary vertices
\cite{dganis}.

Another way to improve the VB limiter is to make it
anisotropic, i.e., to define
$u_{ih}^*$ using the general formula~\eqref{slim} and choose an
individual correction factor $\beta_{ik}$ for each partial
derivative~$U_{ik}$. Such anisotropic limiting approaches were explored, e.g.,
in \cite{dganis,hoteit2004,berger2015}. Their disadvantages (high cost of
solving inequality-constrained optimization problems or the
use of closed-form approximations based on worst-case
assumptions) are similar to those of slope limiting with
nonlinear flux constraints. An additional
disadvantage is the fact that the
 DMP property of cell averages needs to be assumed.

\subsection{Anisotropic slope limiting under derivative constraints}
\label{sec:DClim}

The aforementioned drawbacks of anisotropic slope limiting under solution
constraints of the form~\eqref{vblmp} are due to the fact that the pointwise
value of $u_{ih}^*$ at each vertex $\mathbf{x}_{ip},\ p=1,\ldots,N$ depends on
$d$ correction factors $\beta_{i1},\ldots, \beta_{id}$. Hence, it is
impossible to determine the optimal value of any $\beta_{ik}$ without taking
the other correction factors into account. To avoid this problem, we
introduce a slope limiter which constrains each directional
derivative $U_{ik}^*=\pd{u_{ih}^*}{x_k}=\beta_k\pd{u_{ih}}{x_k},
\ k=1,\ldots,d$ using individually chosen bounds
$U_{ik}^{\max}$ and $U_{ik}^{\min}$.  Instead
of constructing them in a way which
would guarantee the DMP
property for pointwise solution values and/or cell averages, we define
\beq\label{boundsder}
U_{ik}^{\max}=\max_{j\in\mathcal J_{ik}}U^R_{jk},\qquad
U_{ik}^{\min}=\min_{j\in\mathcal J_{ik}}U^R_{jk}
\eeq
using low-order reconstructed values $U^R_{jk}$. By default, the
bounding
stencil $\mathcal J_{ik}$ contains the indices of all common-vertex
neighbors of $K_i$. However, the use of individually chosen
stencils for certain partial derivatives may be
appropriate for anisotropic transport problems
(see Section \ref{sec:num:anis}).

In contrast to finite volume and reconstruction-based DG schemes
in which accurate reconstruction of derivatives from cell
averages is required, a rough approximation $U^R_{jk}$ to
$\pd{u}{x_k}$ in $K_j$ is sufficient to construct reasonable
bounds for anisotropic slope limiting. In view of the fact
that
\beq
\nabla u_{ih}|_{K_i}=\frac{1}{|K_i|}\sum_{j\in\mathcal N_i}
\int_{S_{ij}}u_{ih}\mathbf{n}_{ij}\ds\qquad \forall u_{ih}\in\mathbb{P}_1(K_i)
\eeq
by the divergence theorem, the simple formula
\beq\label{gradrec}
(U^R_{i1},\ldots,U^R_{id})^T=\frac{1}{|K_i|}\sum_{j\in\mathcal N_i}
\int_{S_{ij}}\left(\frac{U_{j0}+U_{i0}}{2}\right)\mathbf{n}_{ij}\ds
\eeq
can be used to calculate the approximate derivatives $U^R_{jk},\ j\in\mathcal J_{ik}$
that define the local bounds \eqref{boundsder}.

The inequality constraints of our new slope limiting procedure are formally defined by
\beq\label{DMPder}
\min\{0,U_{ik}^{\min}\}\le U_{ik}^*=\beta_kU_{ik}\le
\max\{0,U_{ik}^{\max}\}.
\eeq
In practice, we calculate the components $U_{ik}^*$ of the
slope-limited gradient directly as follows:
\beq\label{ansisSL}
U_{ik}^*=\begin{cases}
{\rm minmod}(U_{ik},U_{ik}^{\max}) & \mbox{if}\ U_{ik}>0,\\
0 & \mbox{if}\ U_{ik}=0,\\
{\rm minmod}(U_{ik},U_{ik}^{\min}) & \mbox{if}\ U_{ik}<0.
\end{cases}
\eeq
The  minmod function that we use in this formula is given by
\beq
    {\rm minmod}(a,b)=\begin{cases}
    \min\{a,b\} & \mbox{if}\ a>0,\ b>0,\\
    \max\{a,b\} & \mbox{if}\ a<0,\ b<0,\\
    0 & \mbox{otherwise}.
    \end{cases}
    \eeq

\begin{rmk}    
  The bound-preserving moment limiter proposed by Giuliani and Krivodonova
  \cite{giuliani,giuliani2020} is also based on a comparison of directional
  derivatives to suitably defined reconstructions from cell averages. However,
  our anisotropic slope limiter \eqref{ansisSL} is much simpler, especially in 3D.
  This simplicity is a consequence of the fact that we no longer attempt
  to enforce the DMP constraints \eqref{DMPconstr} for cell averages in the
  process of slope limiting. In our last new method, these constraints are enforced
  using the flux limiters presented in Section \ref{sec:FL}. The purpose
  of slope limiting based on \eqref{DMPder} is to ensure that the difference
  between $U_{ik}^*$ and the inaccurate but consistent reconstruction
  \eqref{gradrec} of the corresponding partial derivative
  from well-behaved BP cell averages cannot become unacceptably large.
\end{rmk}

\begin{rmk} 
The derivative constraints \eqref{DMPder} are far less restrictive than
the solution constraints \eqref{vblmp} of the isotropic VB slope limiter.
For that reason, we do not recommend the use of \eqref{ansisSL} without
flux limiting. The flux limiter should be applied at each SSP Runge-Kutta
stage but slope limiting via \eqref{ansisSL} may be performed just once
per time step or even less frequently. Hence, the combined cost of flux
and slope limiting may, in fact, be lower than that of a more
sophisticated slope limiter.
\end{rmk}

\begin{rmk}
  When it comes to visualization or to computation of derived quantities
  which requires the validity of solution constraints \eqref{vblmp},
  the isotropic VB limiter \eqref{betadef} may be invoked to enforce
  these constraints. The result will be provably BP since the DMP
  property of cell averages is guaranteed by the flux limiter. The
  original values \eqref{ansisSL} of the directional derivatives
  should be used in further computations (if
  any) because the VB postprocessing may introduce additional
  numerical errors.
\end{rmk}

\subsection{Monolithic anisotropic slope limiting}

All slope limiters presented so far were designed to directly adjust the partial
derivatives $U_{ik}$ of a given DG approximation $u_{ih}$. Such predictor-corrector
approaches tend to cause convergence problems at steady state. Since the
solution constraints \eqref{vblmp} and derivative constraints \eqref{DMPder}
do not imply the DMP property of the cell averages, it is not necessary
to enforce these constraints exactly in our method. If the fully discrete
counterpart of \eqref{LLFL} is used as a fixed-point iteration for
solving a stationary hyperbolic problem, the difference between
$U_{ik}$ and $U_{ik}^*$  can be penalized as follows:
\begin{align}
  \int_{K_i}w_{ih}(u_{ih}^{\rm SSP}&-u_{ih})\dx
-\Delta t\int_{K_i}\nabla w_{ih}\cdot\mathbf{f}(u_{ih})\dx\nonumber\\
  &+\Delta t\sum_{j\in\mathcal N_i}\int_{S_{ij}}w_{ih}[\alpha_{ij}
  H(u_{ih},u_{jh},\mathbf{n}_{ij})+(1-\alpha_{ij})
  H_{ij}^{\mathbb{P}_0}]\ds
  \nonumber \\ 
&+\gamma\Delta t
  \left[\int_{K_i}(w_{ih}-W_{i0})(u_{ih}^{\rm SSP} -U_{i0}^{\rm SSP})\dx
   -\int_{K_i}(w_{ih}-W_{i0})(u_{ih}^* -U_{i0})\dx
    \right]
=0,\label{SLviaP}
\end{align}
where $\gamma>0$ is a large penalty parameter, $u_{ih}$ is the solution
at the beginning of the time step or SSP Runge-Kutta stage, $u_{ih}^*$ is
a slope-limited approximation to  $u_{ih}$, and $\alpha_{ij}$
is a flux correction factor that enforces \eqref{DMPconstr} for the
cell average $U_{i0}^*$ of the updated solution $u_{ih}^{\rm SSP}$.
Importantly, the addition of the penalization term does not
change the evolution equation for $U_{i0}^*$ which corresponds
to $w_{ih}\equiv 1$.

\begin{rmk}
  The implicit treatment of the term depending on
  $u_{ih}^{\rm SSP} -U_{i0}^{\rm SSP}$ corresponds to a
  diagonal correction of the element mass matrix
  in DG schemes using the Taylor basis 
  (cf. \cite{entropyDG}).
\end{rmk}  

If flux limiting is performed using formula \eqref{alphaFCT} with the
MCL bounds \eqref{fmaxMCL},\eqref{fminMCL}, the nonlinear
discrete problem defined by \eqref{SLviaP}
has a well-defined residual and steady-state
solutions are independent of the pseudo-time step
$\Delta t$. In Section \ref{sec:num:circ},
we use \eqref{SLviaP} as a fixed-point iteration method.

\section{Entropy stabilization}
\label{sec:ES}

In the case of a nonlinear conservation law, additional corrections
may need to be performed to ensure entropy stability. A convex function
$\eta:\R\to\R$ is called an entropy and $v(u)=\eta'(u)$ is called an
entropy variable if there exists an entropy flux
$\mathbf{q}:\R\to\mathbb{R}^d$ such that
$v(u)\mathbf{f}'(u)=\mathbf{q}'(u)$. A weak solution $u$ of \eqref{goveq}
 is called an entropy solution
if the entropy inequality
\beq
 \pd{\eta}{t}+\nabla\cdot\mathbf{q}(u)\le 0 
\label{ent-ineq}
\eeq
holds for any entropy pair $(\eta,\mathbf{q})$. For any smooth weak
solution, the conservation law
\beq
 \pd{\eta}{t}+\nabla\cdot\mathbf{q}(u)=0 
 \eeq
 can be derived from \eqref{goveq} using multiplication by the entropy
 variable $v$, the chain rule, and the definition of an entropy pair. Hence,
 entropy is conserved in smooth regions and dissipated at shocks.

 A~discretization of \eqref{goveq} is called entropy stable if it satisfies 
 a (semi-)discrete version of the entropy inequality \eqref{ent-ineq}. In the present work, we enforce this property using the
 entropy correction tools developed in \cite{entropyDG}. The
 underlying design criteria impose the following additional
 constraints:
 \begin{itemize}
 \item The limited fluxes $H_{ij}^*$ defined by \eqref{Hstar} should
   satisfy the entropy stability condition
   \beq\label{condES}
   (V_{j0}-V_{i0})H_{ij}^*\le\mathbf{n}_{ij}\cdot
   (\boldsymbol{\psi}(U_{j0})-\boldsymbol{\psi}(U_{i0})),
\eeq
where 
$ \boldsymbol{\psi}(u)=v(u)\mathbf{f}(u)-\mathbf{q}(u)$ and
$V_{i0}=\eta'(U_{i0})$ is an approximation to the entropy variable.

\item The entropy production by solution gradients should be penalized
  using a stabilization term
  \beq
  Q_i(w_{ih},v_{ih})=
  \nu_i\int_{K_i}(w_{ih}-W_{i0},v_{ih}-V_{i0})\dx\ge 0
  \eeq
such that
\beq\label{eprodmax}
P_i(v_{ih},u_{ih})+\sum_{j\in\mathcal N_i}|S_{ij}|G_{ij}^*\le Q_i(v_{ih},v_{ih}),
\eeq
where $P_i(v_{ih},u_{ih})\approx \int_{K_i}\pd{\eta(u)}{t}\dx$
is the rate of entropy production in cell $K_i$ and
\beq
G_{ij}^*=\frac{V_{j0}+V_{i0}}{2}H_{ij}^*-
  \frac12(\boldsymbol{\psi}(U_{j0})
+\boldsymbol{\psi}(U_{i0}))\cdot\mathbf{n}_{ij}
\eeq
is a consistent approximation to the averaged entropy flux
$\frac{1}{|S_{ij}|}\int_{S_{ij}}\mathbf{q}(u)\cdot\mathbf{n}_{ij}\ds$.
 \end{itemize}
 Note that the entropy dissipative term
$Q_i(w_{ih},v_{ih})$ exhibits the same structure as the penalization
terms that we used in \eqref{SLviaP} for monolithic slope
limiting. To avoid severe time step restrictions, we treat
$Q_i(w_{ih},v_{ih})$ implicitly. In the case of a general convex entropy
$\eta(u)\ne\frac{u^2}{2}$, the
need for solving nonlinear systems
is avoided by exploiting the linearized relationship
$$V_{ik}=\eta''(U_{i0})U_{ik},\qquad k=1,\ldots,d$$ between the partial derivatives
of the conserved quantities and entropy variables, see
\cite{entropyDG}.

As shown by Chen and Shu \cite{chen}, the validity of
\eqref{condES} is guaranteed for the low-order LLF flux
$H_{ij}^*=H_{ij}^{\mathbb{P}_0}$. Hence, this condition can
always be satisfied by reducing the value of
$\alpha_{ij}\in[0,1]$ in \eqref{Hstar} if
necessary. In the case $v_{ih}=V_{i0}$, condition
\eqref{eprodmax} holds for $Q_i(v_{ih},v_{ih})=0$. For
$v_{ih}\ne V_{i0}$, the entropy dissipation rate
$Q_i(v_{ih},v_{ih})$ is strictly positive and the
validity of \eqref{eprodmax} can always
be enforced by choosing $\nu_i$ sufficiently large. For
details regarding the definition of entropy stability
preserving flux correction factors $\alpha_{ij}$ and
slope penalization parameters $\nu_i$, we refer the
interested reader to \cite{entropyDG}.

\section{Numerical examples}\label{sec:num}

In this section, we study the numerical behavior of different
slope limiters for DG-$\mathbb{P}_1$ schemes. In our description
of the numerical results, the methods under
investigation are labeled as follows:
\begin{itemize}
\item FC-L: flux-constraining slope limiter \eqref{betaSL1} for linear advection, no flux limiting;
\item FC-N: flux-constraining slope limiter \eqref{betaSL1nonlin} for nonlinear problems, no flux limiting;
\item SC: solution-constraining vertex-based slope limiter \eqref{betadef}, no flux limiting;
\item DC: derivative-constraining anisotropic slope limiter \eqref{ansisSL}
  + MCL flux limiter;
\item DC-M: monolithic version \eqref{SLviaP} of the slope
  limiter \eqref{ansisSL} + MCL flux limiter.
\end{itemize}
A comparative
study of FCT and MCL flux limiters can be found in \cite{EG-MCL}. For
a comparison of the SC approach to other slope limiting
techniques, we refer the reader to Beisiegel \cite{Beisiegel}.

All figures of this section visualize the results obtained on
a mesh of square cells with uniform spacing $h=\frac{1}{128}$.
Numerical solutions are advanced in (pseudo-)time using the
explicit third-order three-stage SSP Runge-Kutta method \cite{ssprev}
and time steps satisfying condition \eqref{cfl}. The
FC-X and SC slope limiters, as well as the MCL flux limiter
and the monolithic version of the DC slope
limiter, are applied after each Runge-Kutta stage. The non-monolithic
DC slope limiter is applied at the beginning of each time step.
Before calculating the error norms and visualizing the
results, we apply the SC limiter. This postprocessing 
constrains the output of all schemes to be in bounds not
only at the centroids but also at the vertices of mesh cells. For
visualization purposes, we project the postprocessed DG-$\mathbb{P}_1$
solutions into the space of continuous bilinear finite elements using the
lumped-mass $L^2$ projection.

\subsection{Anisotropic advection test}\label{sec:num:anis}
\begin{figure}[h!]
\centering

 \small
  
\begin{minipage}[t]{0.5\textwidth}
\centering (a) initial condition \vskip0.2cm

\includegraphics[width=0.85\textwidth,trim=0 20 0 0,clip]{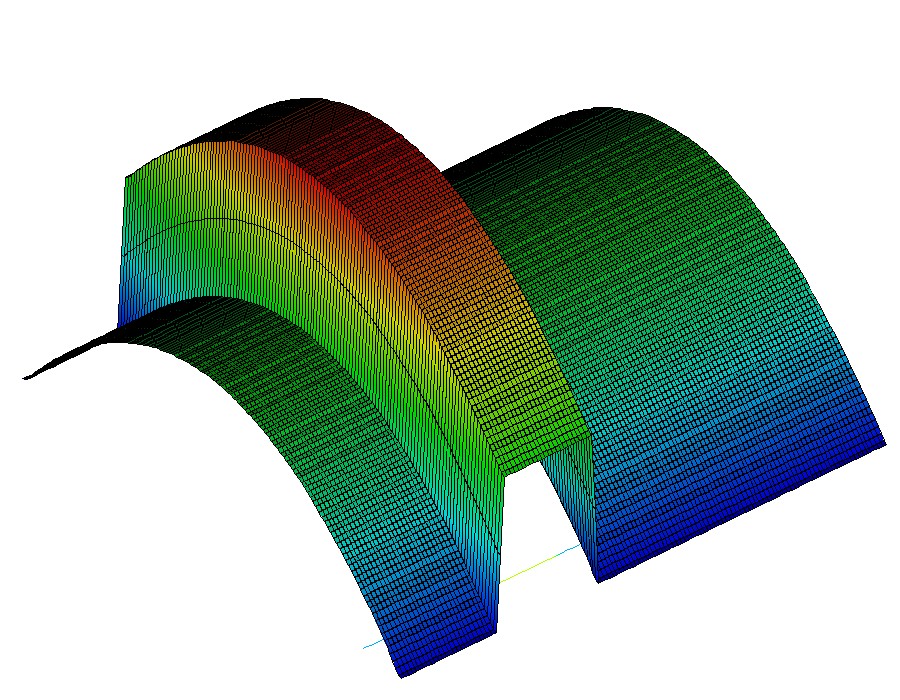}

\end{minipage}%
\begin{minipage}[t]{0.5\textwidth}

  \centering (b) FC-L solution \vskip0.1cm

\includegraphics[width=0.85\textwidth,trim=0 20 0 0,clip]{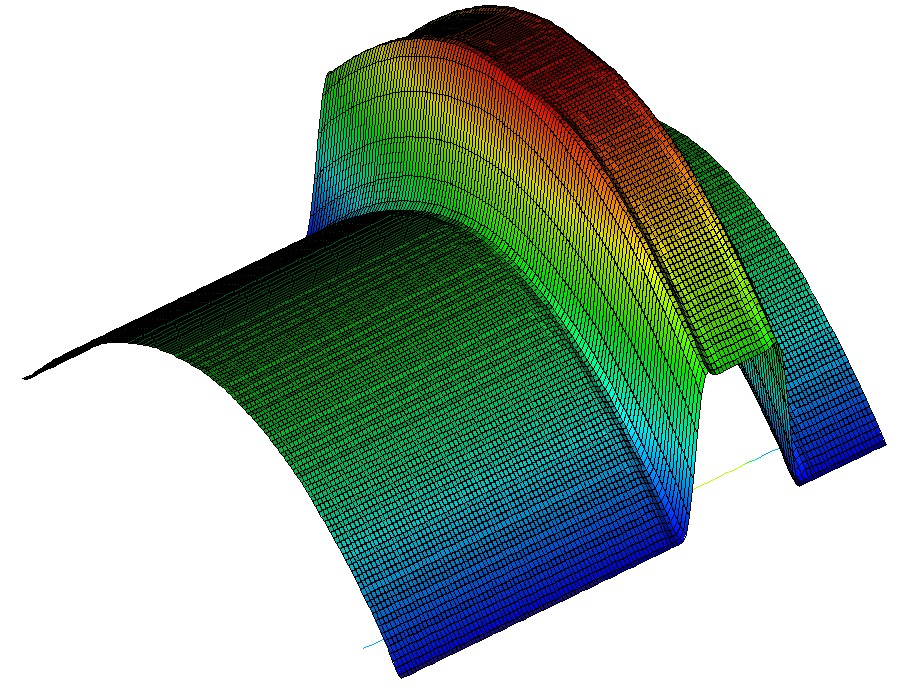}

\end{minipage}

\vskip0.35cm

\begin{minipage}{0.5\textwidth}
\centering (c) SC solution \vskip0.1cm

\includegraphics[width=0.85\textwidth,trim=0 20 0 0,clip]{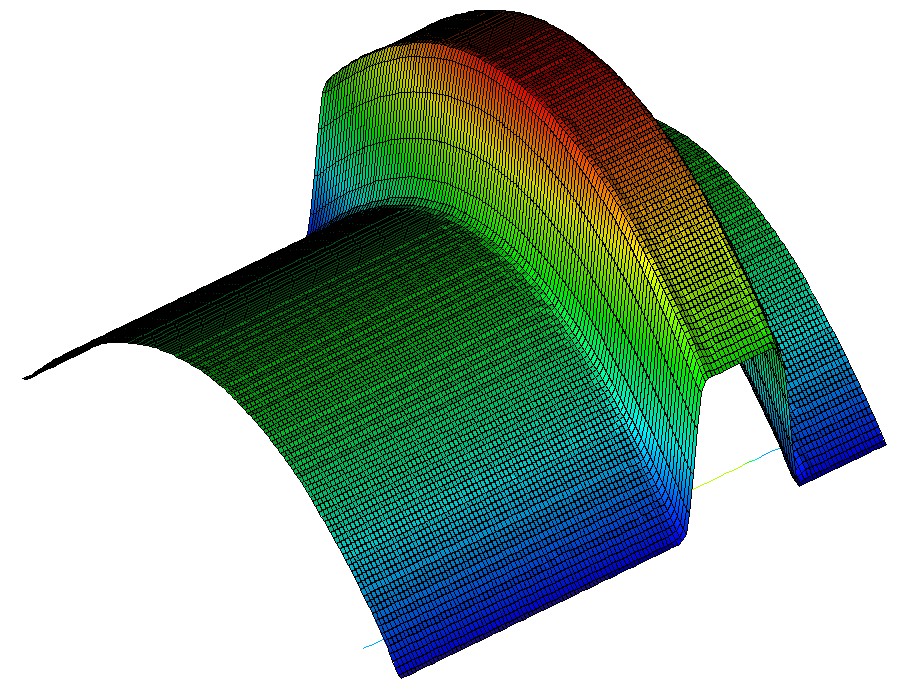}

\end{minipage}%
\begin{minipage}{0.5\textwidth}

\centering (d) DC solution \vskip0.1cm

\includegraphics[width=0.85\textwidth,trim=0 20 0 0,clip]{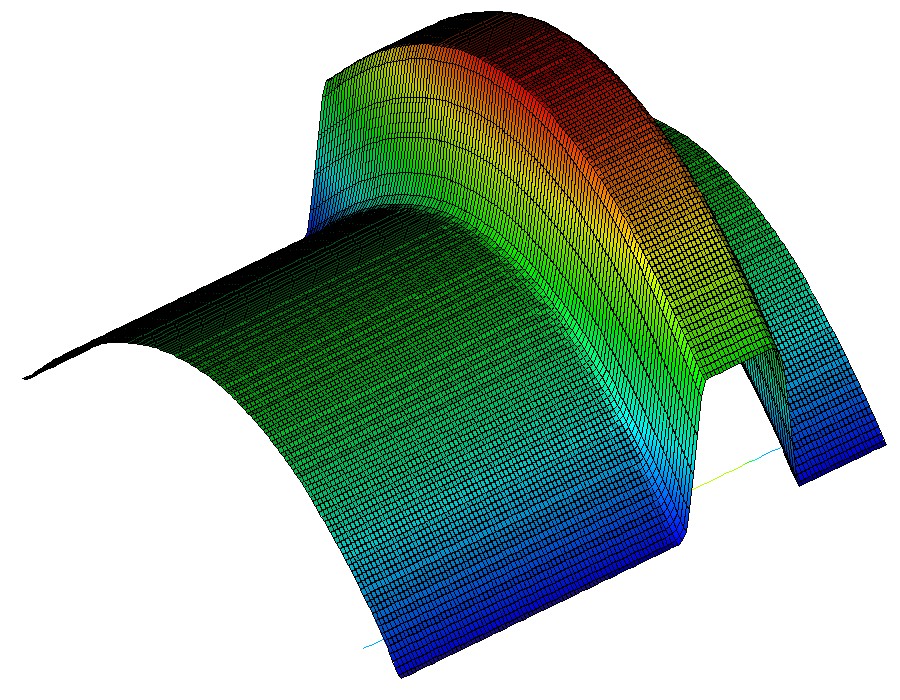}

\end{minipage}

\vskip0.25cm

\caption{Anisotropic advection test. The diagrams show 
  (a) initial data and (b)-(d) numerical solutions at $T=0.4$
  obtained with slope-limited DG-$\mathbb{P}_1$ methods
  using $h=\frac{1}{128}$ and $\Delta t=10^{-3}$.}
\label{AA}

\end{figure}

In the first numerical example, we solve the time-dependent
linear advection equation
\beq
\pd{u}{t}+\nabla\cdot(\mathbf{v}u)=0 \label{linadve}
\eeq
 in $\Omega=(0,1)^2$ using ${\bf v}(x,y)=(0,1)$. The
initial profile displayed in Fig. \ref{AA}(a) is defined by \cite{dganis}
\begin{equation}
u_0(x, y) = w(x) 4 y (1-y), \qquad
\end{equation}
where
\begin{equation}
w(x)=\left\{\begin{array}{ll}
2 & \mbox{if}\ 0.2\le x\le 0.4,\\
1 & \mbox{otherwise}.
\end{array}\right.
\end{equation}
The challenge of this test is to obtain approximations which are
not only DMP-satisfying but also monotonicity-preserving and
do not limit the smooth derivatives $U_{i2}=\pd{u_{ih}}{y}$ in
boundary elements. This example is specifically designed
to show that isotropic limiters may produce inaccurate
results in applications to anisotropic transport problems
such as advection of salinity in ocean flows.

In this experiment, the local bounds \eqref{boundsCA} of
the FC-L slope limiter are defined using cell average data of
all common-vertex neighbors, whereas the bounding stencil
$\mathcal J_{i0}$ of the DC flux limiter is chosen to include
only neighbors belonging to the same horizontal layer of
mesh cells. The bounding stencils $\mathcal J_{i1}$ and
$\mathcal J_{i2}$ of the DC slope limiter for partial
derivatives w.r.t. $x_1:=x$ and $x_2:=y$ are defined using
the horizontal and vertical neighbors, respectively. This
choice is motivated by the anisotropic structure of the
exact solution and exploits the option of anisotropic
limiting in the DC algorithm.

\begin{figure}[h!]
\centering

\includegraphics[width=0.65\textwidth,trim=100 250 100 260,clip]{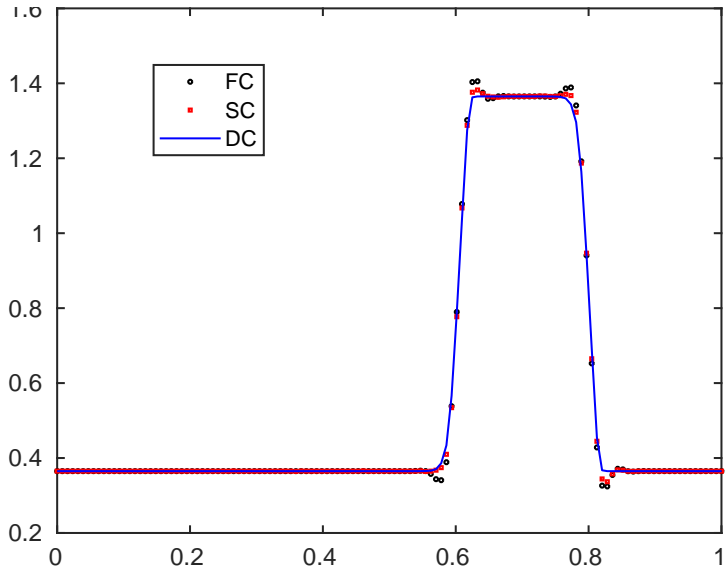}

\caption{Anisotropic advection test. Slope-limited DG-$\mathbb{P}_1$ solutions
  at $T=0.4$ along the line $y=0.1$.}
\label{AA-cut}
 
\end{figure}

Snapshots of slope-limited DG-$\mathbb{P}_1$ solutions
obtained at $T=0.4$ using the time step $\Delta t=10^{-3}$
are presented in Fig.~\ref{AA}. All of them satisfy the
inequality constraints of the underlying limiting
procedures, and the cell averages stay in the range determined
by the maxima and minima of the initial data. However, the
FC-L and SC solutions are not monotonicity-preserving because
they were obtained using multidimensional slope limiters with
isotropic bounds. The local bounds of the FC-L method are
wider than those of the vertex-based SC limiter, which
results in larger undershoots/overshoots around the
discontinuities. The anisotropic bounds of the DC method,
and the way in which it constrains the directional
derivatives, make it possible to preserve monotonicity
along the grid lines. In Fig.~\ref{AA-cut} we compare
the solution profiles produced by the methods under
investigation along the line $y=0.1$. The DC solution
is perfectly monotone because the anisotropic bounds
are tight and the directional derivatives are limited
separately. The isotropic FC-L and SC slope limiters fail
to preserve the 1D structure of the exact solution along
the grid lines because the magnitude of the vertical gradient
affects the bounds for variations of cell averages in the
horizontal direction and vice versa.

\subsection{Solid body rotation}
\label{sec:num:sbr}

In the second linear advection test, we solve \eqref{linadve}
in $\Omega=(0,1)^2$ using the rotating velocity field
${\bf v}(x,y)=(0.5-y,x-0.5)^\top$.
The initial condition, as defined by LeVeque \cite{leveque},
is given by
\begin{equation*}
   u_0(x,y)=\begin{cases}
    u_0^{\rm hump}(x,y)
     &\text{if}\  \sqrt{(x - 0.25)^2 + (y - 0.5)^2}\le 0.15, \\
    u_0^{\rm cone}(x,y)
     &\text{if}\ \sqrt{(x - 0.5)^2 + (y - 0.25)^2}\le 0.15, \\
     1 &\text{if}\ \begin{cases}
       \left(\sqrt{(x - 0.5)^2 + (y - 0.75)^2}\le 0.15 \right) \wedge \\
       \left(|x - 0.5| \ge 0.025 \vee \ y\ge 0.85\right),
     \end{cases}\\
      0 &  \text{otherwise},
    \end{cases}
\end{equation*}
where
\begin{align}
u_0^{\rm hump}(x,y)&=    \frac14 + \frac14 \cos \left(
     \frac{\pi \sqrt{(x - 0.25)^2 + (y - 0.5)^2}}{0.15}\right),\\
u_0^{\rm cone}(x,y)&= 1-\frac{\sqrt{(x - 0.5)^2 + (y - 0.25)^2}}{0.15}.
\end{align}
Homogeneous Dirichlet boundary conditions are prescribed on
portions of $\partial\Omega$ where $\mathbf{v}\cdot\mathbf{n}<0$.

\begin{figure}
  \small
  
\begin{minipage}[t]{0.5\textwidth}
\centering (a) exact solution,\ $u_h\in[0.0,1.0]$\vskip0.2cm

\includegraphics[width=\textwidth,trim=0 20 0 20,clip]{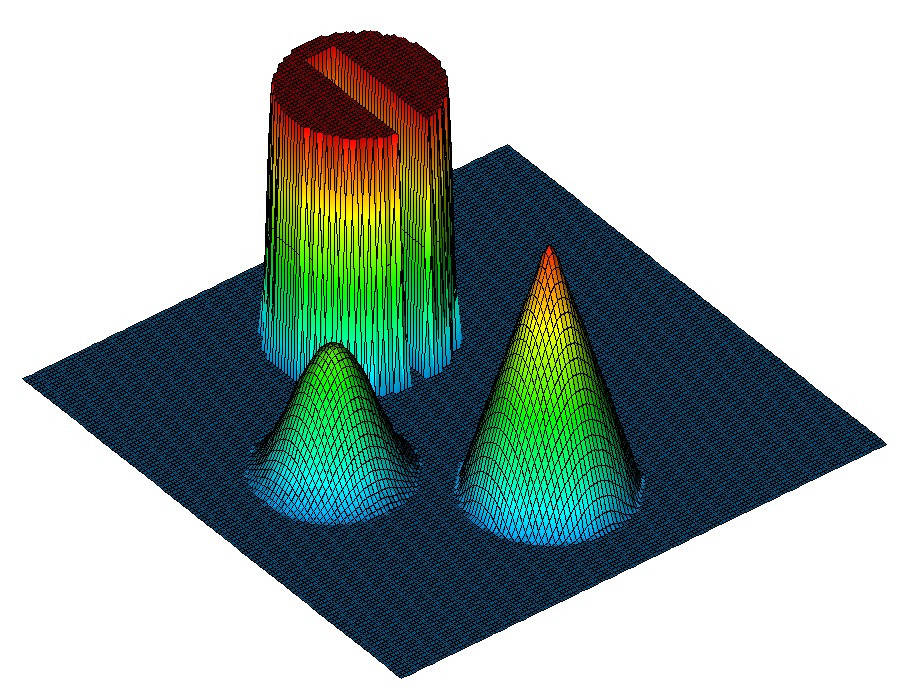}

\end{minipage}%
\begin{minipage}[t]{0.5\textwidth}

  \centering (b) FC-L solution,\ $u_h\in[0.0,0.9995]$\vskip0.2cm

\includegraphics[width=\textwidth,trim=0 20 0 20,clip]{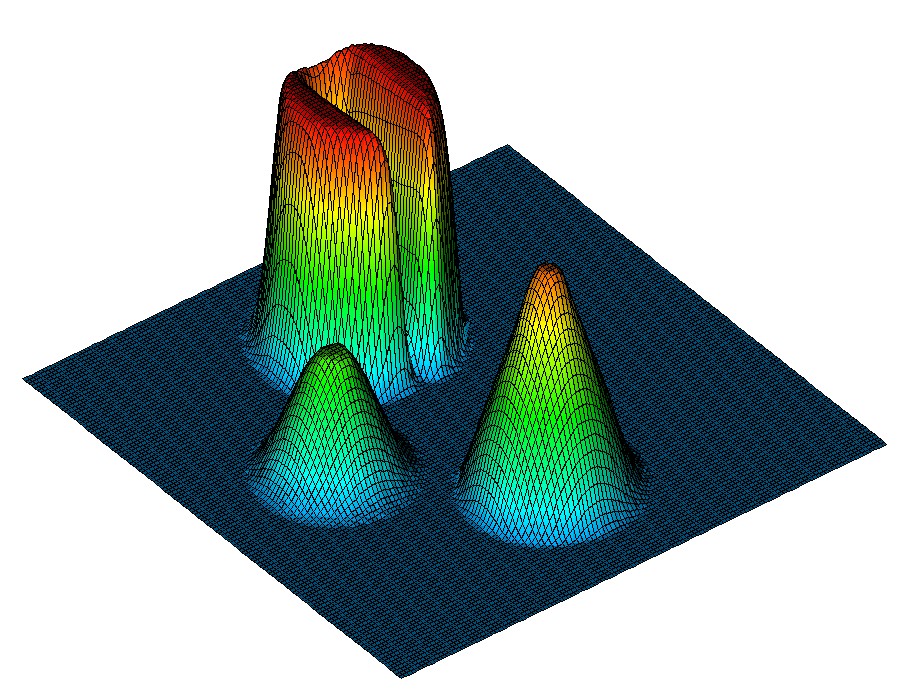}

\end{minipage}

\vskip0.35cm

\begin{minipage}{0.5\textwidth}
\centering (c) SC solution,\ $u_h\in[0.0,0.9956]$\vskip0.2cm

\includegraphics[width=\textwidth,trim=0 20 0 20,clip]{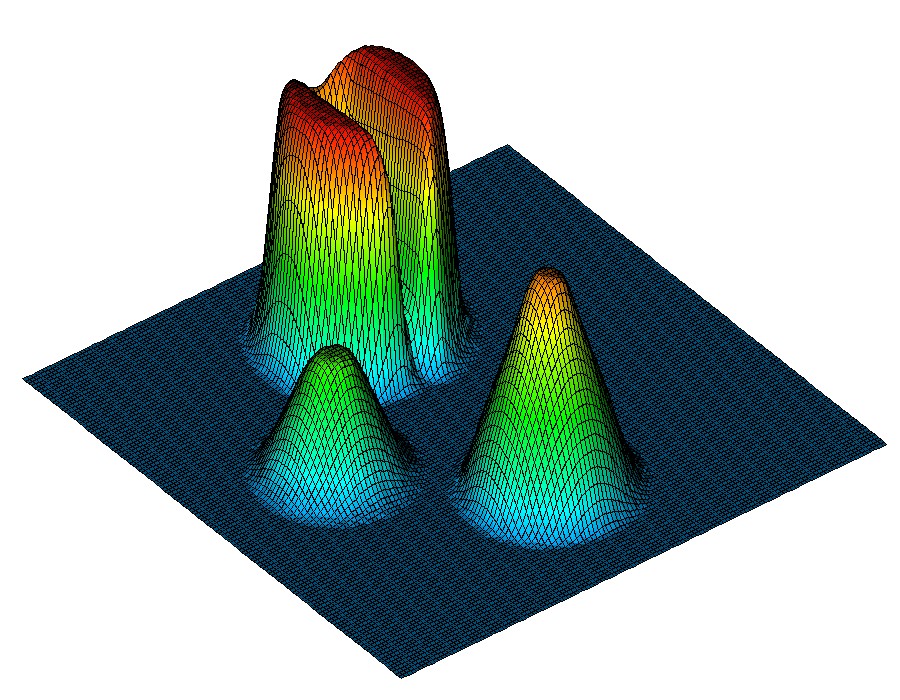}

\end{minipage}%
\begin{minipage}{0.5\textwidth}

\centering (d) DC solution,\ $u_h\in[0.0,0.9998]$\vskip0.2cm

\includegraphics[width=\textwidth,trim=0 20 0 20,clip]{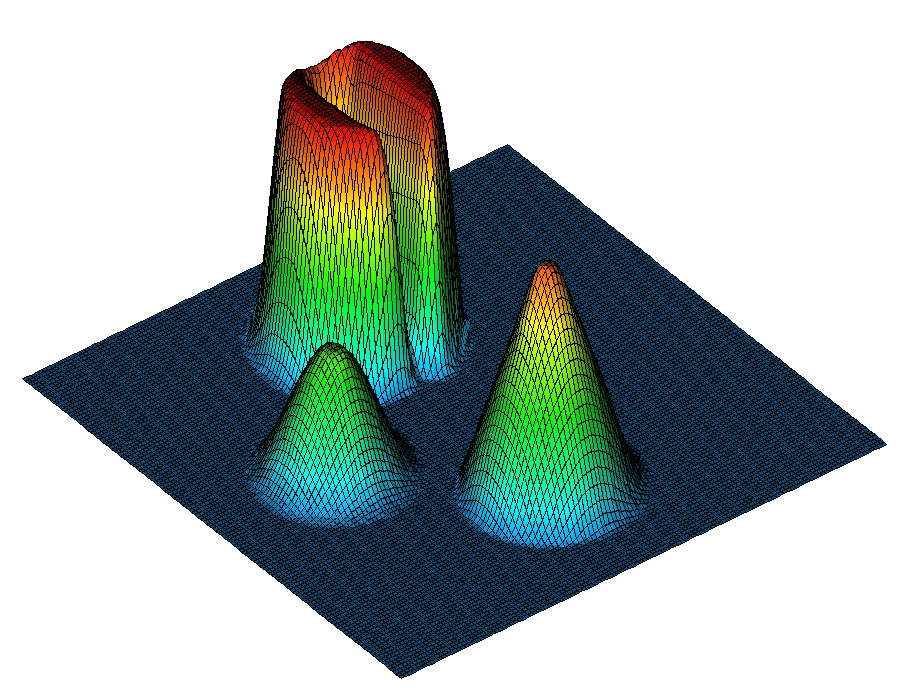}

\end{minipage}

\vskip0.25cm

\caption{Solid body rotation test \cite{leveque}. The diagrams show 
  (a) exact solution / initial data and (b)-(d) numerical solutions
  at $T=2\pi$ obtained with slope-limited DG-$\mathbb{P}_1$ methods
  using $h=\frac{1}{128}$ and $\Delta t=10^{-3}$.}

\label{SBR}

\end{figure}

After each full rotation, the exact solution $u(\cdot,2\pi k),\ k\in
\mathbb{N}$ coincides with the initial data $u_0$. Numerical
solutions are evolved using the time step $\Delta t=10^{-3}$.
In this example and in the remaining test runs of this section,
we use isotropic bounding stencils $\mathcal J_{ik},\ k=0,1,2$
not only for the FC-L limiter but also for the DC version. The
initial data and numerical results produced by the three
methods after one full rotation, i.e., at the final time
$T=2\pi$, are shown in Fig.~\ref{SBR}. The differences
between the numerical solutions are not as pronounced as
in the first example. The
values of the global maxima listed above each plot indicate that
the most diffusive approximation is produced by the vertex-based
SC limiter (which in turn is far less diffusive than many other
limiters for DG schemes, cf. \cite{Beisiegel}). The new FC-L
and DC approaches impose DMP constraints on cell averages
(rather than solution values at the vertices) using local
bounds corresponding to the maximum and minimum of the SC
bounds at the corners of the cell. This limiting strategy
results in better preservation of the initial profile.
In particular, the DC version captures the smooth peak
remarkably well for a method that does not use smoothness
indicators to relax the local bounds. Moreover,
the FC-L and DC methods are provably bound-preserving. The
SC limiter does not produce any undershoots or overshoots
in this particular test but it does not generally guarantee
the BP property as long as the flux limiter is deactivated.

\subsection{Steady circular advection}
\label{sec:num:circ}

In the last linear advection test, we calculate steady-state
solutions to \eqref{linadve} in
$\Omega=(0,1)^2$
using $\mathbf{v}(x,y)=(y,-x)$. The inflow
boundary condition and the exact solution are given by
\begin{equation}
u(x,y)=\left\{\begin{array}{ll}
1, &\quad \mbox{if} \ \ 0.15\le r(x,y)\le 0.45,\\
\cos^2\left(10\pi\frac{r(x,y)-0.7}{3}\right),
&\quad \mbox{if} \ \ 0.55\le r(x,y)\le 0.85,\\
0, & \quad\mbox{otherwise},
\end{array}\right. \label{cc-init}
\end{equation}
where $r(x,y)=\sqrt{x^2+y^2}$ is the distance to the corner point $(0,0)$
of the unit square domain.

\begin{figure}
  \small
  
\begin{minipage}[t]{0.5\textwidth}
\centering (a) exact solution,\ $E_2$=0.000e0

\includegraphics[width=\textwidth,trim=0 20 0 70,clip]{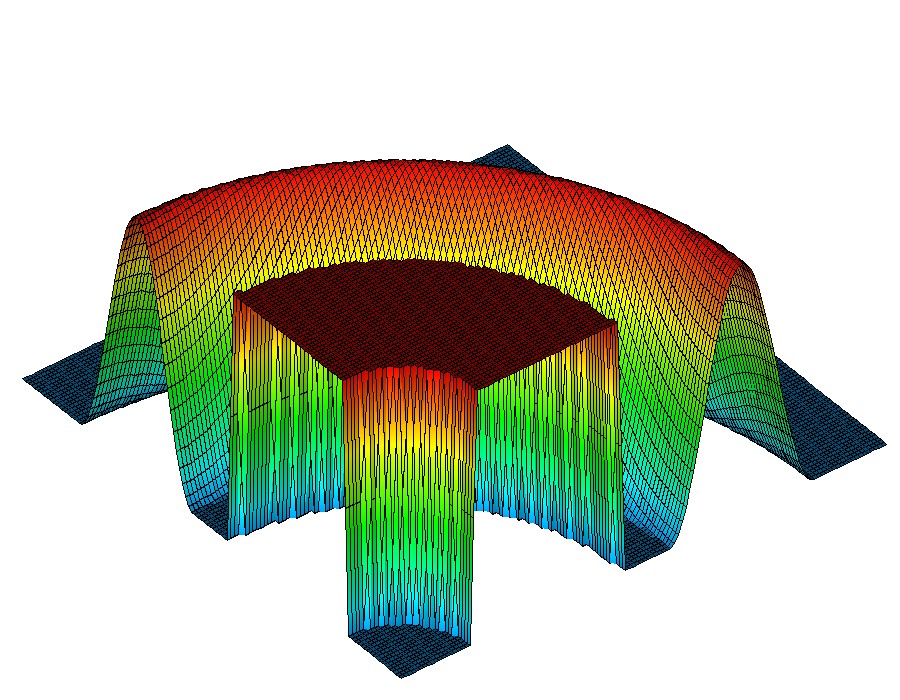}

\end{minipage}%
\begin{minipage}[t]{0.5\textwidth}

  \centering (b) FC-L solution,\ $E_2$=0.397e-01

\includegraphics[width=\textwidth,trim=0 20 0 70,clip]{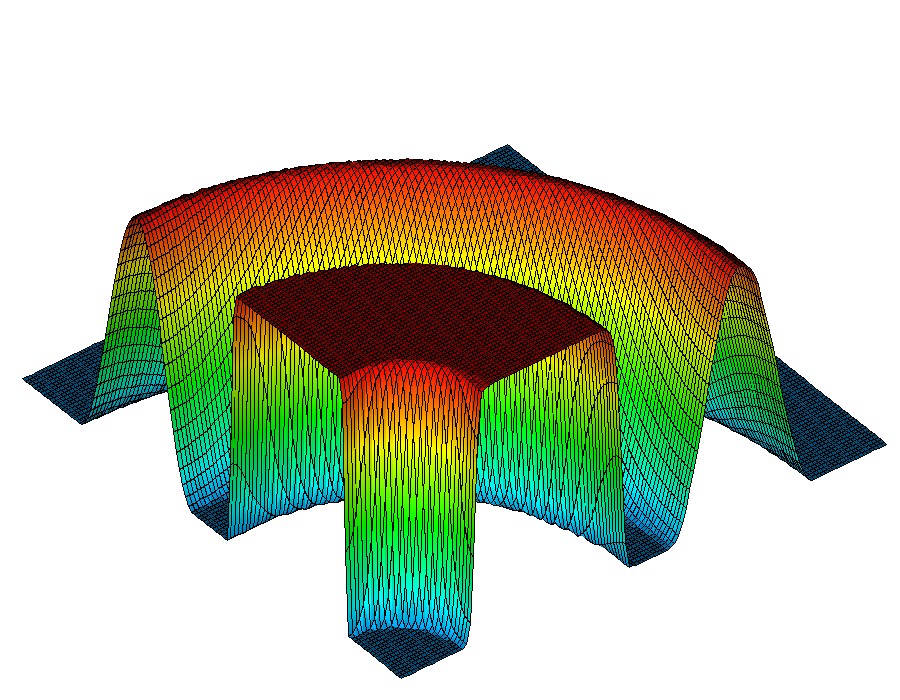}

\end{minipage}

\vskip0.35cm

\begin{minipage}{0.5\textwidth}
\centering (c) SC solution,\ $E_2$=0.415e-01

\includegraphics[width=\textwidth,trim=0 20 0 70,clip]{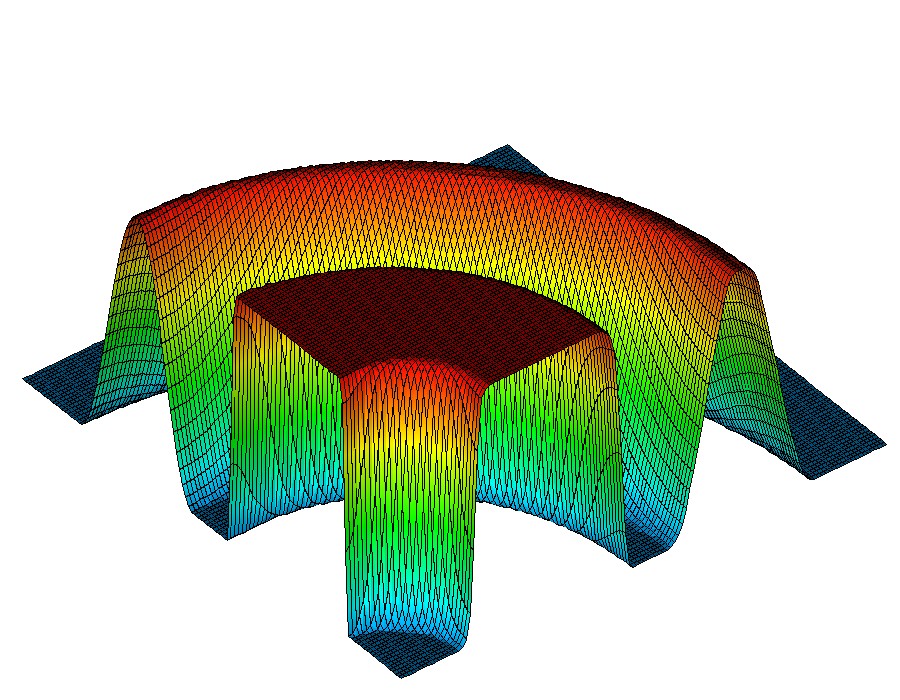}

\end{minipage}%
\begin{minipage}{0.5\textwidth}

\centering (d) DC-M solution,\ $E_2$=0.396e-01

\includegraphics[width=\textwidth,trim=0 20 0 70,clip]{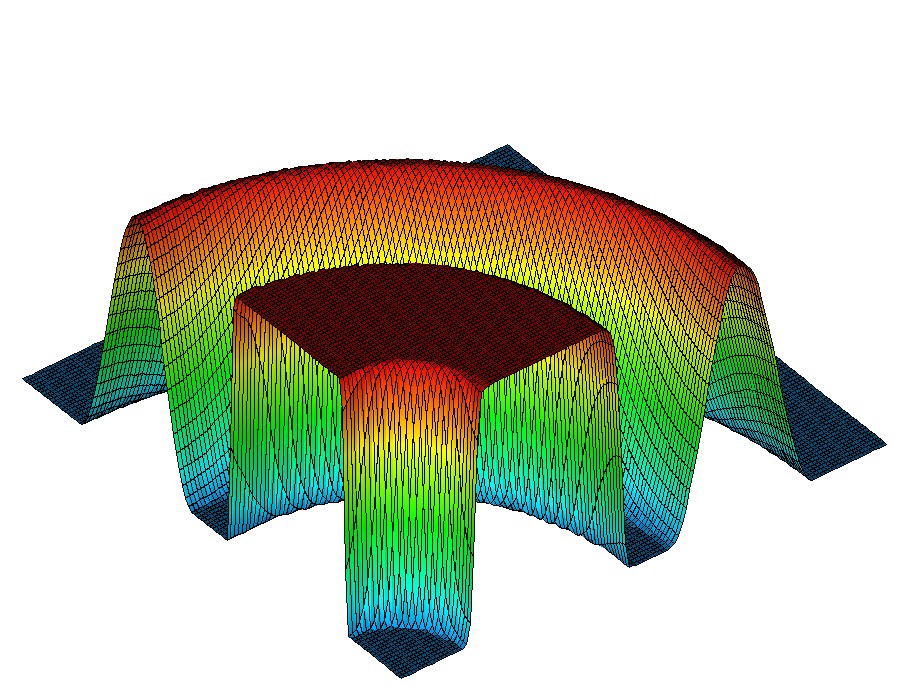}

\end{minipage}

\vskip0.25cm

\caption{Steady circular advection test \cite{convex,EG-MCL}. The diagrams show 
  (a) exact solution and (b)-(d) numerical solutions obtained with slope-limited
  DG-$\mathbb{P}_1$ methods using $h=\frac{1}{128}$.}

\label{CC}

\end{figure}

Numerical solutions are marched to the steady state with the same explicit SSP Runge-Kutta scheme that we use for time-dependent problems. The methods under investigation are FC-L, SC, and DC-M. The parameter $\gamma$ for the penalty term of the monolithic DC scheme \eqref{SLviaP} is chosen to be $10^{3}\cdot M_i$, where $M_i$ is a diagonal entry of the Taylor mass matrix (see \cite{dglim,entropyDG} for the definition of the Taylor basis). In our experience,  the DC-M results are quite insensitive to the choice of $\gamma$ because the main bound-preserving correction tool of this method is the MCL flux limiting procedure.  
 
Computations are terminated when the $L^1$ norm of the difference $U_{i0}^*-U_{i0}$ becomes smaller than the prescribed tolerance $10^{-10}$ at the first SSP RK stage of the pseudo-time step. We remark that we were often unable to reach this tolerance with the SC approach. It is widely known that SC-like limiters may inhibit convergence to steady-state solutions. The monolithic FC-L and DC-M schemes have well-defined residuals and, therefore, exhibit better steady state convergence behavior.

The exact solution \eqref{cc-init}
and  numerical solutions obtained with the three methods are displayed in Fig.~\ref{CC}. For a better quantitative comparison, we list the $L^2$ errors in the approximation of cell averages (denoted by $E_2$) above each plot. The FC-L and DC-M errors are virtually the same, while the SC error is slightly larger due to more restrictive inequality constraints. All solutions are bounded by the maximum $u^{\max}=1$
and minimum $u^{\min}=0$
of the Dirichlet boundary data. No spurious oscillations are observed in the neighborhood of discontinuities, and smooth portions of the exact solution are reproduced very well.
To study the numerical behavior of the three limiting techniques for a
steady advection test with a globally smooth exact solution, we perform grid
convergence studies for 
\cite{CL-diss}
\beq\label{expIC}
u(x,y)=\exp\left(-100(r(x,y)-0.7)^2\right),\qquad 0\le x,y\le 1
\eeq
using the same velocity field. In Table~\ref{CCconv}, we present the
$E_2$ convergence history for the standard DG-$\mathbb{P}_0$ method and
its slope-limited counterparts. The experimental order of convergence
(EOC) is above
2.0 if no limiting is performed and approximately 1.75 otherwise.
Further improvements can be achieved by using less restrictive bounds,
as proposed in \cite{timelim}. The remarkably similar convergence
behavior of the FC-L and DC schemes and the larger absolute errors of 
the SC method indicate that the choice of the local bounds has a stronger impact
on the accuracy of slope-limited approximations to smooth solutions
than the way in which these bounds are enforced in our methods. 

\begin{table}[h!]
\centering

\tabcolsep0.4cm
\renewcommand{\arraystretch}{1.1}

\begin{tabular}{|c|cc|cc|cc|cc|}
\hline
$h$ & $E_2^{\mathbb{P}_1}$ & EOC
& $E_2^{\rm FC}$ & EOC &  $E_2^{\rm SC}$ & EOC &
$E_2^{\rm DC}$ & EOC
\\\hline
$\frac{1}{32}$  & 7.17e-3 &      & 1.39e-2 &      &  3.28e-2 &      & 1.38e-2 &      \\
$\frac{1}{64}$  & 1.68e-3 & 2.09 & 4.02e-3 & 1.79 &  1.04e-2 & 1.66 & 4.00e-2 & 1.79 \\
$\frac{1}{128}$ & 4.09e-4 & 2.04 & 1.20e-3 & 1.75 &  3.00e-3 & 1.79 & 1.20e-3 & 1.74 \\
$\frac{1}{256}$ & 1.01e-4 & 2.02 & 3.59e-4 & 1.74 &  8.30e-4 & 1.85 & 3.59e-4 & 1.74  
\\\hline
\end{tabular}

\vskip0.25cm

\caption{ Convergence
  behavior of the DG-$\mathbb{P}_1$ methods in steady-state
  computations for the steady circular advection test with the
  smooth inflow profile \eqref{expIC}.}
\label{CCconv}

\end{table}

\subsection{KPP problem}
The KPP problem \cite{Guermond2016,Guermond2017,kpp} is a
challenging nonlinear test for verification of entropy
stability properties. In this last example, we solve
\eqref{goveq} in $\Omega=(-2,2)\times(-2.5,1.5)$
using the flux function
\beq
\mathbf{f}(u)=(\sin(u),\cos(u))
\eeq
and the initial
condition
\beq
u_0(x,y)=\begin{cases}
\frac{14\pi}{4} & \mbox{if}\quad \sqrt{x^2+y^2}\le 1,\\
\frac{\pi}{4} & \mbox{otherwise}.
\end{cases}
\eeq
The entropy flux corresponding to $\eta(u)=\frac{u^2}2$ is
$\mathbf{q}(u)=(u\sin(u)+\cos(u),u\cos(u)-\sin(u))$.
The exact solution exhibits
a two-dimensional rotating wave structure which is difficult to
capture correctly.

\begin{figure}
  \small
  
\begin{minipage}[t]{0.5\textwidth}
  \centering (a) DG-$\mathbb{P}_0$ solution

\includegraphics[width=0.9\textwidth,trim=100 20 100 20,clip]{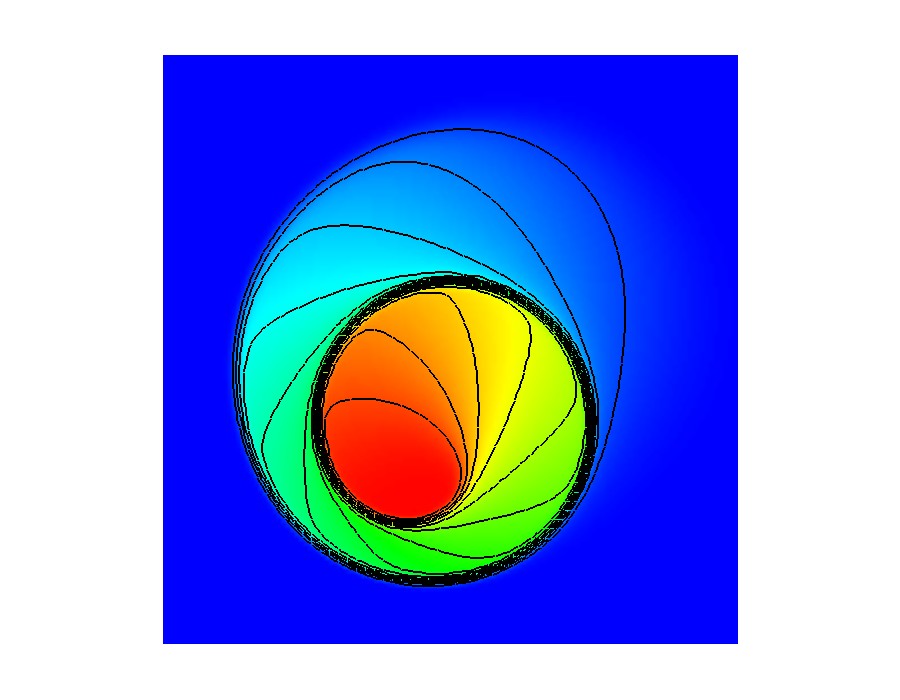}

\end{minipage}%
\begin{minipage}[t]{0.5\textwidth}

\centering (b) FC-N solution

\includegraphics[width=0.9\textwidth,trim=100 20 100 20,clip]{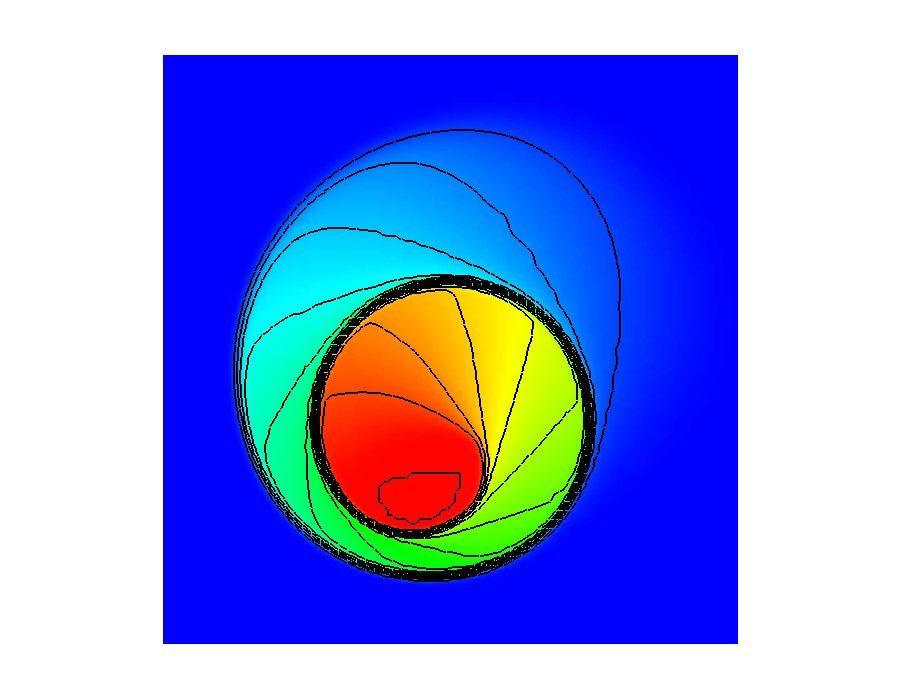}

\end{minipage}

\begin{minipage}[t]{0.5\textwidth}

  \centering (c) SC solution

\includegraphics[width=0.9\textwidth,trim=100 20 100 20,clip]{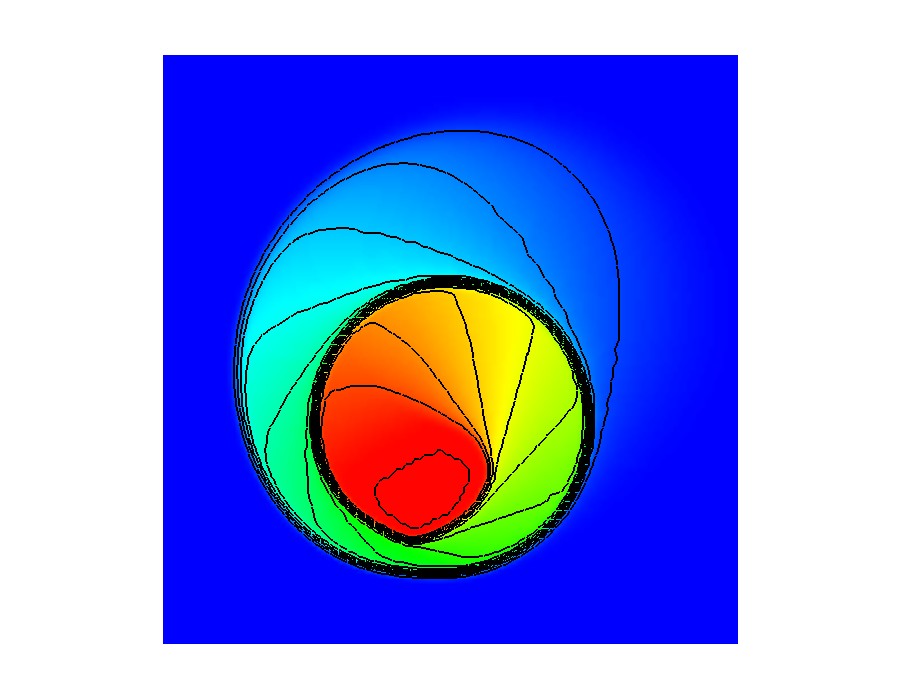}

\end{minipage}%
\begin{minipage}[t]{0.5\textwidth}

  \centering (d) DC solution

\includegraphics[width=0.9\textwidth,trim=100 20 100 20,clip]{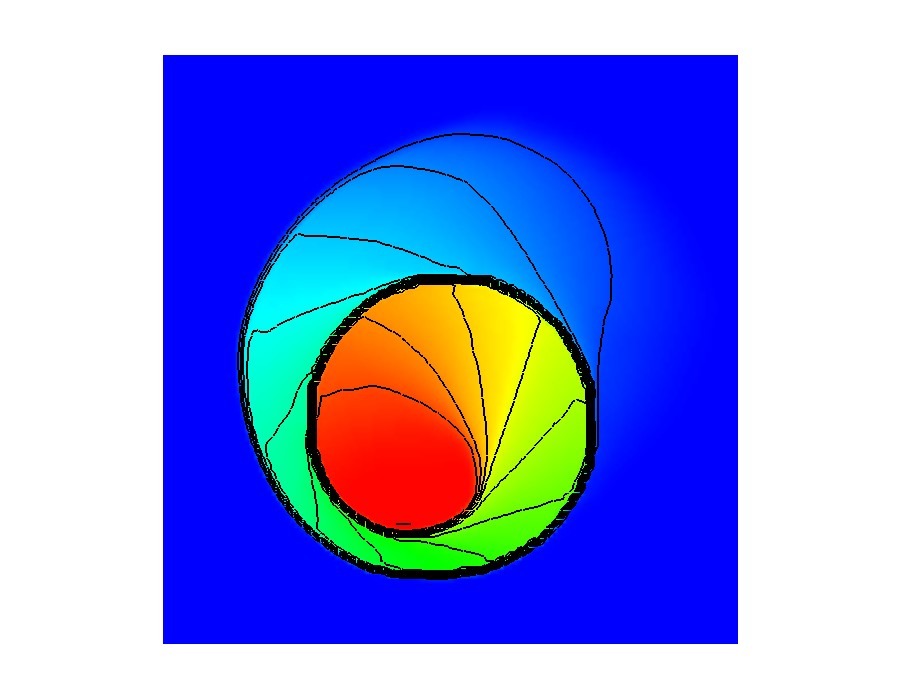}

\end{minipage}

\vskip0.25cm

\caption{KPP problem \cite{kpp}. The 2D plots show (a)
  property-preserving
  DG-$\mathbb{P}_0$ solution and (b)-(d) entropy-stabilized
  slope-limited
  DG-$\mathbb{P}_1$ solutions at $t=1.0$ obtained with
  $h=\frac{1}{128}$ and $\Delta t=10^{-3}$.}

\label{KPP2D}

\end{figure}

\begin{figure}[h!]
  \small
  
\begin{minipage}[t]{0.5\textwidth}
\centering (a) DG-$\mathbf{P}_0$ solution \vskip0.2cm

\includegraphics[width=\textwidth,trim=0 20 0 20,clip]{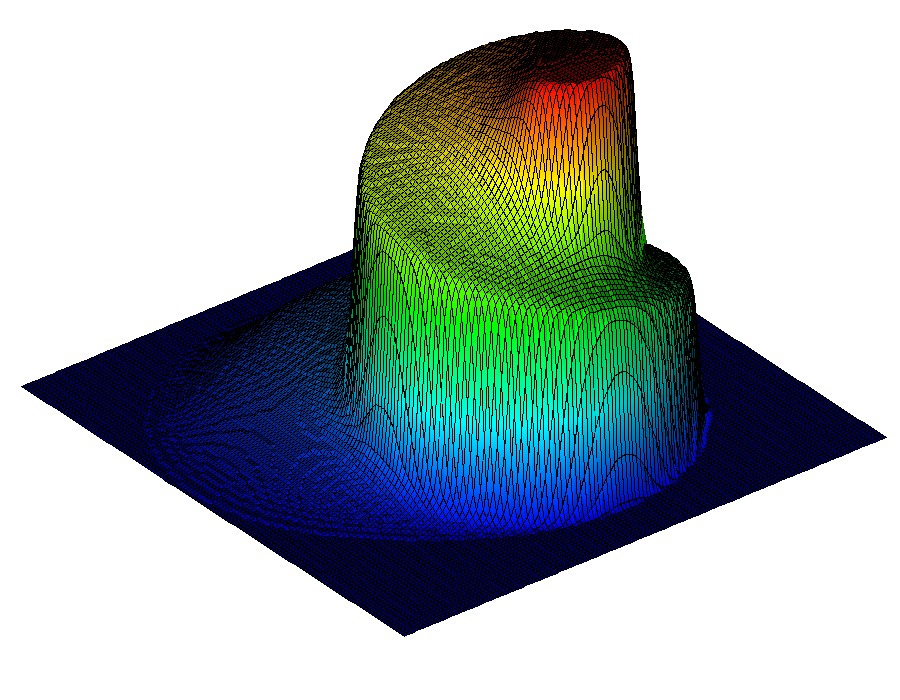}

\end{minipage}%
\begin{minipage}[t]{0.5\textwidth}

\centering (b) FC-N solution \vskip0.2cm

\includegraphics[width=\textwidth,trim=0 20 0 20,clip]{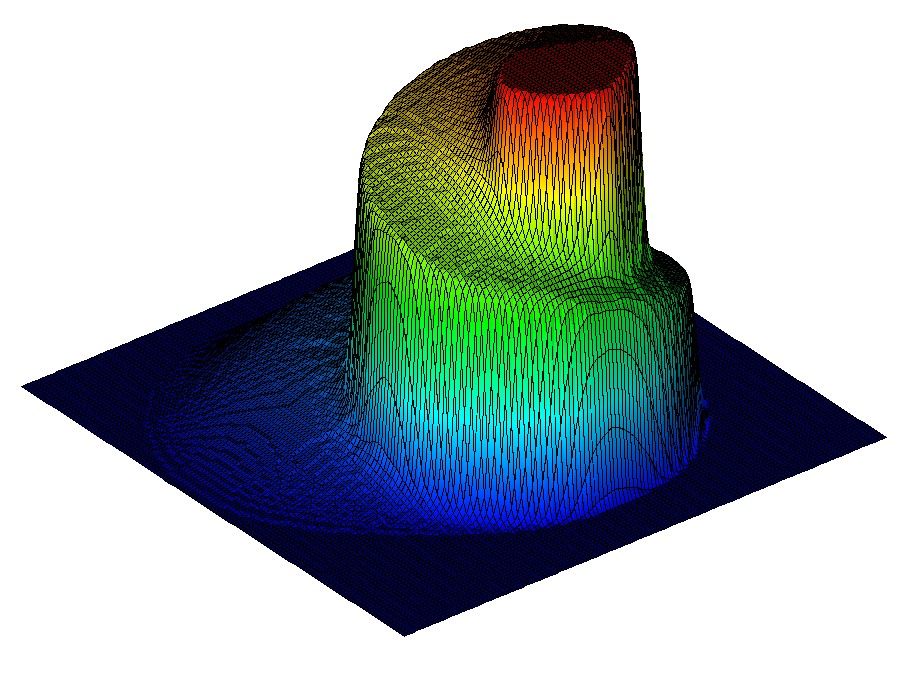}

\end{minipage}

\begin{minipage}[t]{0.5\textwidth}

\centering (c) SC solution \vskip0.2cm

\includegraphics[width=\textwidth,trim=0 20 0 20,clip]{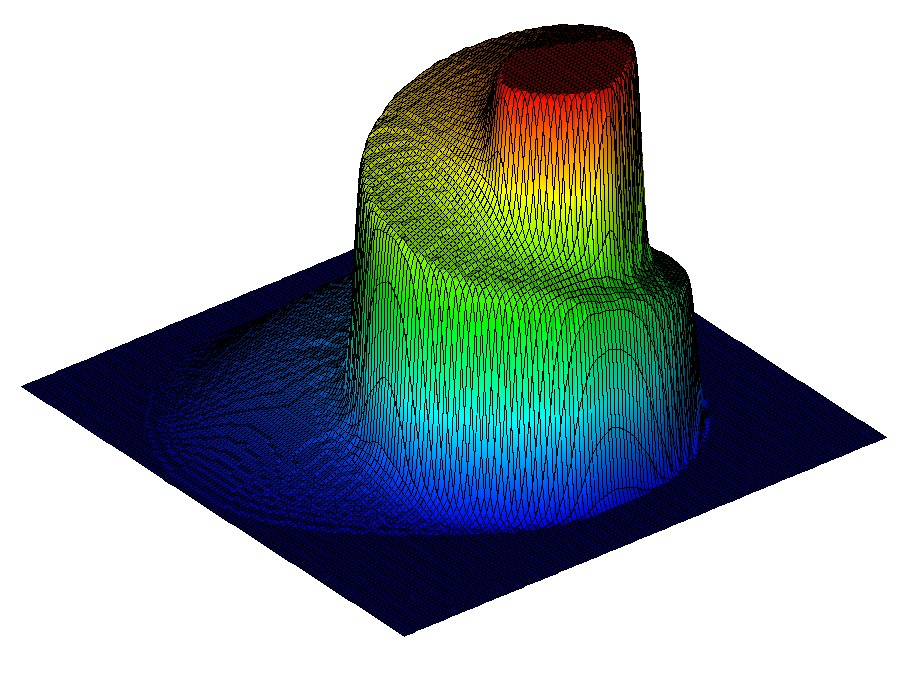}

\end{minipage}%
\begin{minipage}[t]{0.5\textwidth}

\centering (d) DC solution \vskip0.2cm

\includegraphics[width=\textwidth,trim=0 20 0 20,clip]{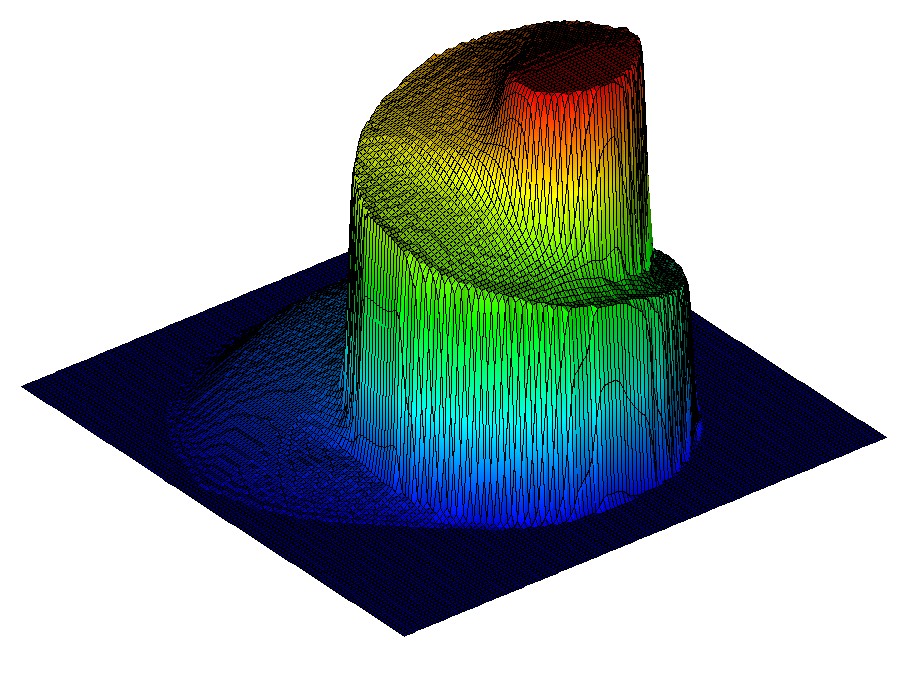}

\end{minipage}

\vskip0.25cm

\caption{KPP problem \cite{kpp}. The 3D plots show (a) property-preserving
  DG-$\mathbb{P}_0$ solution and (b)-(d) entropy-stabilized
  slope-limited DG-$\mathbb{P}_1$ solutions at $T=1.0$ obtained with
  $h=\frac{1}{128}$ and $\Delta t=10^{-3}$.}

\label{KPP3D}

\end{figure}

We define the LLF flux \eqref{LLF_flux} using the global
upper bound $\lambda^{\max}=1$ for the maximum wave speed.
More accurate estimates
can be found in \cite{Guermond2017}. 
In our numerical study of the FC-N, SC, and DC
methods for the KPP problem,
we use the entropy correction tools outlined in
Section~\ref{sec:ES} in addition to bound-preserving
flux/slope limiting. Snapshots of the numerical solutions
at the final time $T=1.0$ are presented
in Figs.~\ref{KPP2D} and \ref{KPP3D}. The diffusive DG-$\mathbb{P}_0$
results are included to illustrate the correct wave structure of the
entropy solution. As noticed by Guermond et al. \cite{Guermond2016},
high-order finite element schemes may converge to wrong weak solutions
of the KPP problem even if they are equipped with bound-preserving
limiters. None of our entropy-stabilized DG-$\mathbb{P}_1$ methods
exhibits this behavior. The twisted shocks stay clearly separated, while
the levels of numerical dissipation are reduced compared to DG-$\mathbb{P}_0$.
The sharpest resolution of discontinuities is obtained with our DC
approach. The new FC-N limiter performs better than the SC limiter.
These encouraging results demonstrate that the DMP property and
entropy stability of a slope-limited DG-$\mathbb{P}_1$ approximation
can be guaranteed using fairly simple limiting techniques which are
at least as robust and accurate as existing alternatives.

\section{Conclusions}
\label{sec:conc}

The findings reported in this work reveal some interesting relationships between flux and slope limiting in DG-$\mathbb{P}_1$ methods for hyperbolic problems. In particular, it turns out that a carefully designed slope limiter can act as a flux limiter, i.e., provably enforce flux constraints which imply a discrete maximum principle for cell averages. However, our preferred limiting strategy is a combination of flux correction for cell averages and slope correction for directional derivatives. We have shown that this approach makes it possible to cure some alarming deficiencies of existing limiting techniques. We also discussed the design of anisotropic and monolithic limiters in this framework. Last but not least, we addressed the aspects of entropy stabilization via flux limiting and slope penalization.

The FCT and MCL flux limiters presented in this work have already been extended to nonlinear hyperbolic systems and guarantee preservation of invariant domains \cite{Guermond2018,Guermond2019,HennesDG,convex,pazner}. Extensions to high-order finite element approximations and general high-order Runge-Kutta methods are feasible as well. In principle, the methodology proposed in the present paper is directly applicable to finite elements of arbitrary order. However, recent advances in the development of algebraic limiting approaches \cite{HennesDG,Hajduk2020,convex2,entropyHO,hpfem,CG-BFCT,pazner} indicate that a localization to subcells is required to achieve at least the same accuracy as with the DG-$\mathbb{P}_1$ scheme using the same number of degrees of freedom. The simplest way to meet this requirement is to construct an $hp$-adaptive partition of unity for the pair of finite element spaces corresponding to a high-order DG method and the DG-$\mathbb{P}_1$ subcell approximation on a submesh with the same nodes \cite{hpfem} . Using a smoothness indicator to select the appropriate local basis in each cell, the application of limiters can be restricted to $\mathbb{P}_1$ subcells without losing the high accuracy of the DG approximation elsewhere. To preserve the high accuracy and DMP property of the space discretization, time integration may need to be performed using a flux-corrected Runge-Kutta method of sufficiently high order. The first representatives of such methods were recently developed in \cite{timelim,CG-BFCT}. In summary, the proposed methodology can be extended to high-order space-time discretizations but many additional aspects must be taken into account to reap the potential benefits.



\begin{thebibliography}{99}

\bibitem{dganis}
 V. Aizinger, A. Kosik, D. Kuzmin, and B. Reuter, Anisotropic slope
 limiting for discontinuous Galerkin methods. {\it Int. J. Numer.
  Methods Fluids} {\bf 84} (2017) 543-565.
  
\bibitem{barthjesp}
T.~Barth and D.C.~Jespersen, 
The design and application of upwind schemes on unstructured meshes.
{\it AIAA Paper}, 89-0366, 1989.


\bibitem{barthohl}
T.~Barth and M. Ohlberger, Finite volume methods: foundation and 
analysis. In: E.~Stein, R.~de~Borst, T.J.R.~Hughes (eds),
{\it Encyclopedia of Computational Mechanics, 
Volume 1: Fundamentals}. John Wiley \& Sons, 2004, 439--474.

\bibitem{Beisiegel}
  N. Beisiegel, {\it High-order Adaptive Discontinuous Galerkin Inundation Modeling}.
  PhD thesis, University of Hamburg, 2014.

\bibitem{shasta}
J.P.~Boris and D.L.~Book, Flux-Corrected Transport: I. SHASTA, a 
fluid transport algorithm that works. {\it J. Comput. Phys.} 
{\bf 11} (1973) 38--69.


\bibitem{BurchardRennau2008} H. Burchard and H. Rennau, Comparative 
quantification of physically and numerically induced mixing in ocean models. 
{\it Ocean Modelling} {\bf 20} (2008) 293--311.


\bibitem{chen}
  T. Chen and C.W. Shu, Entropy stable high order discontinuous Galerkin methods with suitable quadrature rules for hyperbolic conservation laws.
  {\it J. Comput. Phys.} {\bf 345} (2017) 427--461.

\bibitem{christov2008new}
  I. Christov and B. Popov, New non-oscillatory central schemes on unstructured
  triangulations for hyperbolic systems of conservation laws.
  {\it J. Comput. Phys.}
  {\bf 227-11},
  (2008) 5736--5757.


\bibitem{RKDG-II}
B. Cockburn and C.-W. Shu, TVB Runge-Kutta local projection discontinuous Galerkin finite element
method for scalar conservation laws II: General framework. {\em Math. Comp.} {\bf 52} (1989) 411--435.


\bibitem{cotter}
  C.J. Cotter and D. Kuzmin, 
Embedded discontinuous Galerkin transport schemes with localised limiters.
 {\it J. Comput. Phys.} {\bf 311}  (2016) 363--373.

 \bibitem{CH-Paper}
   F. Frank, A. Rupp, and D. Kuzmin, Bound-preserving flux limiting schemes
   for DG discretizations of conservation laws with applications to the Cahn–Hilliard equation
{\it Computer Methods Appl. Mech. Engrg.} {\bf 359} (2020) 112665.

  \bibitem{giuliani}
    A. Giuliani and L. Krivodonova. A moment limiter for the discontinuous Galerkin method on unstructured triangular meshes. {\it SIAM J. Sci. Comput.} {\bf 41} (2019) A508--A537.

      \bibitem{giuliani2020}
    A. Giuliani and L. Krivodonova. A moment limiter for the discontinuous Galerkin method on unstructured tetrahedral meshes. {\it J. Comput. Phys.} {\bf 404} (2020) 109106.
 
\bibitem{ssprev}
S.~Gottlieb, C.-W.~Shu, and E.~Tadmor, Strong stability-preserving
high-order time discretization methods. {\it SIAM Review}
{\bf 43} (2001) 89--112.

\bibitem{Griffies2000} 
 S.M. Griffies, R.C. Pacanowski, and R. W. Hallberg,
 Spurious diapycnal mixing associated with advection in a
 $z$-coordinate ocean model. 
{\it Monthly Weather Review} {\bf 128} (2000) 538--564. 


\bibitem{Guermond2018}
J.-L. Guermond, M. Nazarov, B. Popov, and I. Tomas, Second-order
invariant domain preserving approximation of the Euler equations
using convex limiting. {\it SIAM J. Sci. Computing} {\bf 40}
(2018) A3211-A3239.

\bibitem{Guermond2019}
  J.-L. Guermond, M. Nazarov, and I. Tomas, Invariant domain preserving
  discretization-independent schemes and convex limiting for hyperbolic
  systems. {\it Computer Methods Appl. Mech. Engrg.} {\bf 347} (2019)
  143--175.

\bibitem{Guermond2016}
  J.-L. Guermond and B.~Popov, Invariant domains and first-order continuous
  finite element approximation for hyperbolic systems.
   {\it SIAM J. Numer. Anal.} {\bf 54} (2016) 2466--2489.

 \bibitem{Guermond2017}
     J.-L. Guermond and B.~Popov, 
     Invariant domains and second-order continuous finite element approximation for scalar conservation equations. {\it SIAM J. Numer. Anal.} {\bf 55} (2017)
     3120--3146.

\bibitem{HennesDG}
H. Hajduk, Monolithic convex limiting in discontinuous Galerkin discretizations of hyperbolic conservation laws.
Preprint {\tt  arXiv:2007.01212v2 [math.NA]}, 2020.

\bibitem{Hajduk2020}
  H. Hajduk, D. Kuzmin, Tz. Kolev, V. Tomov, I. Tomas, and J.N. Shadid,
  Matrix-free subcell residual distribution for Bernstein finite
  elements: Monolithic limiting.
  {\it Computers \& Fluids}. {\bf 200} (2020) 104--451.

\bibitem{harten1}
A.~Harten, High resolution schemes for hyperbolic conservation laws.
{\it J. Comput. Phys.}  {\bf 49} (1983) 357--393.

\bibitem{harten2}
A.~Harten, On a class of high resolution total-variation-stable
finite-difference-schemes. {\it SIAM J. Numer. Anal.} {\bf 21} (1984) 1-23.
  
\bibitem{hoteit2004}
H. Hoteit, Ph. Ackerer, R. Mos\'e, J. Erhel, and B. Philippe, 
New two-dimensional slope limiters for discontinuous Galerkin 
methods on arbitrary meshes. {\it Int. J. Numer. Meth. Engrg.}
{\bf 61} (2004) 2566--2593. 


\bibitem{jameson1}
A.~Jameson, Computational algorithms for aerodynamic analysis and
design. {\it Appl. Numer. Math.} {\bf 13} (1993) 383--422.

\bibitem{jameson2}
A.~Jameson, Analysis and design of numerical schemes for gas dynamics
1. Artificial diffusion, upwind biasing, limiters and their effect on
accuracy and multigrid convergence. {\it Int. Journal of CFD} 
{\bf 4} (1995) 171--218.

\bibitem{kriv2007}
L.~Krivodonova, Limiters for high-order discontinuous Galerkin 
methods. {\it J. Comput. Phys.} {\bf 226} (2007) 879--896.
  
\bibitem{kriv2004}
L.~Krivodonova, J.~Xin, J.-F.~Remacle, N.~Chevaugeon, and 
J.E.~Flaherty, Shock detection and limiting with discontinuous 
Galerkin methods for hyperbolic conservation laws. 
{\it Appl.~Numer.~Math.} {\bf 48} (2004) 323--338. 
  
\bibitem{kpp}
A. Kurganov, G. Petrova, and B. Popov, Adaptive semidiscrete central-upwind
schemes for nonconvex hyperbolic conservation laws. {\it SIAM J. Sci.
Comput.} {\bf 29} (2007) 2381--2401.


 \bibitem{dglim}
    D. Kuzmin, A vertex-based hierarchical slope limiter for 
    p-adaptive discontinuous Galerkin methods. {\it J. Comput. 
    Appl. Math.} {\bf 233} (2010) 3077--3085.

 \bibitem{entropyDG} 
   D. Kuzmin, Entropy stabilization and property-preserving limiters for discontinuous Galerkin discretizations of nonlinear hyperbolic equations. Preprint
  {\tt  arXiv:2004.03521 [math.NA]}, 2020.
  
\bibitem{convex}
  D. Kuzmin,  Monolithic convex limiting for continuous finite element
  discretizations of hyperbolic conservation laws. {\it Comput.
    Methods Appl. Mech. Engrg.} {\bf 361} (2020) 112804.

 \bibitem{EG-MCL}
   D. Kuzmin, H. Hajduk, and A. Rupp, Locally bound-preserving
   enriched Galerkin methods for the linear advection equation.
{\it Computers and Fluids}
 {\bf 205} (2020) 104525.

 \bibitem{convex2} 
   D. Kuzmin and M. Quezada de Luna, Subcell flux limiting for high-order Bernstein finite element discretizations of hyperbolic conservation laws.
   {\it J. Comput. Phys.} {\bf 411} (2020) 109411.

    \bibitem{entropyHO} 
      D. Kuzmin and M. Quezada de Luna, Entropy conservation property and entropy stabilization of high-order continuous Galerkin approximations to scalar conservation laws. Preprint {\tt arXiv:2005.08788 [math.NA]}, 2020.
      
\bibitem{hpfem}
 D. Kuzmin, M. Quezada~de Luna, and C. Kees, A partition of unity approach to
adaptivity and limiting in continuous finite element methods.
{\it Computers \& Mathematics with Applications} {\bf  78}
(2019) 944--957.

 \bibitem{timelim}
   D. Kuzmin, M. Quezada~de Luna, D. Ketcheson, and J. Gr\"ull,
   Bound-preserving convex limiting for high-order 
    Runge-Kutta time discretizations of hyperbolic
  conservation laws. In preparation.
   
\bibitem{leveque92}
R.J.~LeVeque, {\it Numerical Methods for Conservation
Laws}. Birkh\"auser, 1992.


\bibitem{leveque}
R.J. LeVeque, High-resolution conservative algorithms for advection in incompressible flow.
{\it SIAM Journal on Numerical Analysis} {\bf 33}, (1996) 627--665.


\bibitem{CL-diss}
C. Lohmann, {\it Physics-Compatible Finite Element Methods
  for Scalar and Tensorial Advection Problems.}
Springer Spektrum, 2019.

\bibitem{FCTsync}
C. Lohmann and D. Kuzmin, Synchronized flux limiting for gas dynamics
variables. {\it J. Comput. Phys.} {\bf  326} (2016) 973--990.

\bibitem{CG-BFCT}
 C. Lohmann, D. Kuzmin, J.N. Shadid, and S. Mabuza, Flux-corrected transport
algorithms for continuous Galerkin methods based on high order Bernstein
finite elements. {\it J. Comput. Phys.} {\bf 344} (2017) 151-186.

\bibitem{berger2015}
S. May and M. Berger, Two-dimensional slope limiters for 
finite volume schemes on non-coordinate-aligned meshes.  
{\it SIAM J. Sci. Comput.} {\bf 35} (2013)  A2163--A2187. 

   \bibitem{Moe2017}
S.A. Moe, J.A. Rossmanith, and D.C. Seal, 
Positivity-preserving discontinuous Galerkin methods with 
Lax-Wendroff time discretizations.  {\it J. Sci. 
Comput.} {\bf 71} (2017) 44-70.

\bibitem{pazner}
W. Pazner, Sparse invariant domain preserving discontinuous Galerkin methods with subcell convexlimiting. Preprint {\tt arXiv:2004.08503 [math.NA]}, 2020.

  \bibitem{zalesak79}
S.T.~Zalesak, Fully multidimensional flux-corrected transport
algorithms for fluids. {\it J. Comput. Phys.} {\bf 31} (1979) 335--362.


\bibitem{zalesak87}
S.T.~Zalesak, A preliminary comparison of modern shock-capturing
schemes: linear advection. In: R.~Vichnevetsky and R.~Stepleman
(eds), {\it Advances in Computer Methods for PDEs}. 
Publ. IMACS, 1987, 15--22.

\bibitem{zhang1}
X. Zhang and C-W. Shu, On positivity-preserving high order discontinuous 
Galerkin schemes for compressible Euler equations on rectangular meshes.
{\it J. Comput. Phys.} {\bf 229} (2010) 8918--8934.

\bibitem{zhang2}
X. Zhang and C-W. Shu, Maximum-principle-satisfying and
positivity-preserving high-order schemes for conservation laws: survey
and new developments. {\it Proc. R. Soc. A} {\bf 467} (2011) 2752--2776.
  
\end{thebibliography}
\end{document}